\documentclass[leqno]{article}
\usepackage{latexsym,amsmath,amsfonts,mathrsfs,amssymb,amsthm}
\usepackage[usenames]{color}
\usepackage{graphicx}
\def\R{\mathbb R}
\def\N{\mathbb N}

\def\Z{\mathbb Z}

\def\a{\alpha}
\def\e{\epsilon}
\def\d{\delta}

\def\be{\begin{equation}}
\def\ee{\end{equation}}
\def\bs{\backslash}
\def\qed{\hfill$\Box$\bigskip}
\def\nd{\noindent Proof. }
\numberwithin{equation}{section}
\newtheorem{lem}[equation]{Lemma}
\newtheorem{pro}[equation]{Proposition}
\newtheorem{defn}[equation]{Definition}
\newtheorem{thm}[equation]{Theorem}
\newtheorem{cor}[equation]{Corollary}
\newtheorem{rem}[equation]{Remark}
\begin{document}
\bigskip

\centerline{\Large \textbf{Topological minimal sets and their applications}}

\bigskip

\centerline{\large Xiangyu Liang}

\vskip 1cm

\centerline {\large\textbf{Abstract.}}

\bigskip

In this article we introduce a definition of topological minimal sets, which is a generalization of that of Mumford-Shah-minimal sets. We prove some general properties as well as two existence theorems for topological minimal sets. As an application we prove the topological minimality of the union of two almost orthogonal planes in $\R^4$, and use it to improve the angle criterion under which the union of several  higher dimensional planes is Almgren-minimal.

\bigskip


\bigskip


\section{Introduction}

\bigskip

One of the canonical topics in the geometric measure theory is the theory of minimal sets. Briefly a minimal set is just a set which minimizes the Hausdorff measure among a certain class of competitors. Different choice of class of competitors gives different kind of minimal sets. So we have the following general definition. 

\begin{defn}[Minimal sets]\label{min}
Let $0<d<n$ be integers. A closed set $E$ in $\R^n$ is said to be minimal of dimension $d$ in $\R^n$ if 
\be H^d(E\cap B)<\infty\mbox{ for every compact ball }B\subset U,\ee
and
\be H^d(E\bs F)\le H^d(F\bs E)\ee
for every competitor $F$ for $E$.
\end{defn}

Among all class of minimal sets, the notion of Almgren minimal sets, which has been introduced by F.Almgren \cite{AL76}, leads to important progress in understanding the regularity of minimal surfaces and existence results for Plateau's problem (especially, the classification of singularities of soap films, due to Jean Taylor, 1976, which says that every Al-minimal set of dimension 2 in $\R^3$ is locally $C^1$ equivalent to a minimal cone. c.f. \cite{Ta}). So let us recall here the global version of the definition of Almgren competitors.

\begin{defn}[Almgren competitor (Al-competitor for short)] Let $E$ be a closed set in $\R^n$. An Almgren competitor for $E$ is a closed set $F\subset \R^n$ that can be written as $F=\varphi(E)$, where $\varphi:\R^n\to\R^n$ is a Lipschitz map such that there exists a compact ball $B\subset\R^n$ such that
\be \varphi_{B^C}=id\mbox{ and }\varphi(B)\subset B.\ee 
Such a $\varphi$ is called a deformation in $B$, and $F$ is also called a deformation of $E$ in $B$.
\end{defn}

An Almgren minimal set is a set as in Definition \ref{min} where we take Al-competitors as the class of competitors. Intuitively, a $d-$dimensional Al-minimal set is a closed set whose $d-$Hausdorff measure could not be decreased by any local Lipschitz deformation.

This definition is intuitive enough to describe the behavior of soap films and is easy to understand. However, the definition of Al-competitors depends too much on the parametrization (because it uses deformations), which, sometimes brings non-necessary obstacles to the verifications of some nice properties of Al-minimizers. 

Another class of minimal sets----Mumford-Shah (MS) minimal sets, is less intuitive. Let us give first the definition of Mumford-Shah competitors.

\begin{defn}[Mumford-Shah (MS) competitors]Let $E$ be a closed set in $\R^n$. A Mumford-Shah competitor for $E$ is a set $F\subset\R^n$ such that there exists a compact ball $B\subset \R^n$ verifying

1) $F\bs B=E\bs B$;

2) For all $y,z\in\R^n\bs(B\cup E)$ who are separated by $E$, $y,z$ are also separated by $F$.

Here ``$y,z$ are separated by $E$" means that $y$ and $z$ belong to two different connected components of $\R^n\bs E$.
\end{defn} 

We use the name "Mumford-Shah" because the separation condition 2) comes from the study of minimal segmentations in the Mumford-Shah functional (c.f.\cite{DMS,DJT}). This condition is an essentially topological condition which depends only on the set itself. A Mumford-Shah minimal set is just a minimal set in Definition 1.1 where we take MS competitors as the class of competitors. But here, the separation condition implies that the set $E$ has to be of codimension 1 (otherwise no separation occurs). 

This topological condition intervenes indirectly in the story of soap films. We can prove that for a fixed set $E$, deformations keeps the separation condition (c.f.\cite{Du} XVII 4.3), so all its Al-competitors are MS-competitors. Therefore all MS-minimizers are Al-minimizers. We do not know yet whether the two classes of minimizers are equal. But at least in the whole $\R^3$, the answer is very likely to be yes. (For further information, see \cite{DJT}).

An obvious advantage of defining a competitor by posing some topological condition like separation is that an existence theorem for a minimizer is easier to get. In fact we only have to verify that the topological condition is stable under passing to the limit. This is much easier than proving that a limit of deformations is still a deformation, since the later needs some intelligent parametrization.

MS-minimizers admit another good property that we are not yet able to prove for Al-minimizers. That is, the product of a MS-minimizer (in $\R^n$) with $\R^m$ is still MS-minimal in $\R^{n+m}$. (This is an example where deformations bring some trouble for Al-minimizers.) In fact an intuitive efficient way of constructing new minimal sets is by taking the product of two minimal sets. But the minimality of the product of two minimal sets is still an open question.  

As we have pointed out, the notion of MS-minimizers exists only in codimension 1. So in the rest of this article we will try to generalize this notion to higher codimensions  properly, so as to keep all its good properties.  The seperation condition coincides with two topological invariants (which are also defined for higher codimension): homotopy groups and homology groups in codimension 1.  We prefer taking homology groups, that is

\begin{defn}[Topological competitors] Let $E$ be a closed set in $\R^n$. We say that a closed set $F$ is a topological competitor of dimension $d$ ($d<n$) of $E$, if there exists a ball $B\subset\R^n$ such that

1) $F\bs B=E\bs B$;

2) For all Euclidean $n-d-1$-sphere $S\subset\R^n\bs(B\cup E)$, if $S$ represents a non-zero element in the singular homology group $H_{n-d-1}(\R^n\bs E;\Z)$, then it is also non-zero in $H_{n-d-1}(\R^n\bs F;\Z)$.
\end{defn}

And Definition 1.1 gives the definition of topological minimizers.

The reason for which we don't take homotopy groups instead to define topological competitors is the following. We wish that the class of topological minimizers is a subclass of Al-minimizers, as is the class of MS-minimizer. The easiest way is to prove that all Al-competitors are topological competitors. Remark 3.6 tells us that this is not true if we use homotopy groups.  But luckily we have the chance with the homology groups. We will prove in Section 3 that every Almgren competitor $F$ of a closed set $E$ is a topological competitor of $E$ (c.f. Proposition 3.7).

Thus we decide to use homology groups. 

We will also prove in Section 3.23 that the product of a topological minimal set of dimension $d$ with $\R^m$ is a topological minimal set of dimension $d+m$ (c.f. Proposition 3.23).

\medskip

In Section 4 we will prove two theorems of existence of topological minimal sets. We cite one of them here, since they are kind of similar.

\noindent{\textbf{Theorem 4.2.}} {\it  Let $\Omega\subset\R^n$ be a compact set that admits a Lipschitz neighborhood retraction from $\Omega_\e$ to $\Omega$ with $\e>0$. Let $h$ be a continuous function from $\Omega$ to $\R$, and suppose that $1\le h\le M$ on $\Omega$. Let $\{w_j\}_{j\in J}$ be a family of smooth $n-d-1-$ surfaces in $\R^n\bs\Omega_\e$ which are non-zero in $H_{n-d-1}(\R^n\bs \Omega)$. Set
 \be \mathfrak F=\{F\subset\Omega\mbox{ closed, and for all }j\in J, w_j\mbox{ is non-zero in }H_{n-d-1}(\R^n\bs F)\}.\ee
Define
\be J_h(F)=\int_Fh(x)dH^d(x).\ee
Then there exists $F_0\in\mathfrak F$ such that
\be J_h(F_0)=\inf\{J_h(F);F\in\mathfrak F\}.\ee}

In addition, in Section 5 we will give an application of this generalization, which is related to a conjecture given by F.Morgan \cite{Mo93}. He gave an conjectural angle criterion under which the union of two $m-$planes is Al-minimal, for $m\ge 2$. We recall (part of) this conjecture here.

If $P,Q$ are two $m-$dimensional planes in $\R^{2m}$, we can describe their relative geometric position by their characteristic angles $(\a_1,\cdots,\a_m)$, with $0\le\a_1\le\cdots\le\a_m\le\frac\pi2$: among all pairs of unit vectors $v\in P$ and $w\in Q$, choose and fix $v_1,w_1$ which minimize the angle between them. Next we choose $v_2\in P,w_2\in Q$ with $v_2\perp v_1,w_2\perp w_1$, which minimize the angle among all such pairs. We repeat this $m$ times. For each $i$, denote by $\a_i$ the angle between $v_i$ and $w_i$. These are the characteristic angles between $P$ and $Q$.

Morgan's angle conjecture is the following.
The union of two $m-$planes $P_1\cup P_2$, with characteristic angles $\a_1,\cdots,\a_m$, is Almgren minimal if and only if the two conditions below are satisfied:

1) $\alpha_m\le\alpha_1+\cdots+\alpha_{m-1}+\frac\pi3$ ; 2) $\alpha_1+\cdots+\alpha_m\ge\frac{2\pi}{3}.$

In fact 2) is a necessary condition, as proved by Gary Lawler \cite{La94}. The remaining part is still open. 

 In Section 5 we will prove that the union of two almost orthogonal planes in $\R^4$ is topologically minimal. Then as a corollary of this and Proposition 3.23, the product of $\R^n$ with a union of two almost orthogonal planes is also topologically minimal, and hence Al-minimal. Recall that in \cite{XY10} we have proved that the union of two almost orthogonal (and hence transversal) $m-$planes is Al-minimal. So this last result of the paper tells us that when $m>2$, two $m-$planes do not need to be transversal when their union is Al-minimal. This makes a progress in the direction of the conjecture. However between almost orthogonal and non transversal, there is a large gap, that we still do not understand what happens.
 
 \textbf{Acknowledgement:} I would like to thank Guy David for many helpful discussions and for his continual encouragement. I also wish to thank Pierre Vogel for having given many useful suggestions. The results in this paper were part of the author's doctoral dissertation at University of Paris-Sud 11 (Orsay). This work was supported by grants from R\'egion Ile-de-France.

\section{Definitions and notations on topology}

In this section we give some necessary definitions and notations of topology. 

\subsection{Algebraic topology}

\subsubsection{Simplicial homology (with coefficients in $\Z$)}

\begin{defn}[Simplex] An $n$-simplex $s$ is the convex hull of $n+1$ ordered points $v_0\cdots, v_n$ which can be used to form an affine coordinate system in an Euclidean space of dimension $n$. The convex hull of any subset (ordered according to their order in the simplex) of these $n+1$ points is called a face of this simplex. Those faces are simplices themselves. In particular, the convex hull of a subset of size $m+1$ (of these $n+1$ points) is a $m$-simplex, called a $m-$face of the initial simplex.
\end{defn}

\begin{defn}[Simplicial complex] A simplicial complex $K$ is a set of simplices such that

$1^\circ$ All faces of each simplex in $K$ belong to $K$;

$2^\circ$ Any non-empty intersection of any two simplices $\sigma_1,\sigma_2\in K$ is a face of both $\sigma_1$ and $\sigma_2$.

A $k-$simplicial complex $K$ is a simplicial complex where the largest dimension of any simplex in $K$ is $k$. We say a complex $K$ is finite if $K$ is a finite set.
\end{defn}

\begin{defn}[Polyhedron and triangulation] The polyhedron associated to a simplicial complex $K$, denoted $|K|$, is the union of all simplices belonging to $K$. A pair $(K,\pi)$ of a simplicial complex $K$ and a homeomorphism between its polyhedron $|K|$ and a topological space $X$ is called a triangulation of $X$. We say that the triangulation is finite if $K$ is finite.
\end{defn}

\begin{defn}[Simplicial chain] Let $K$ be a simplicial complex. A simplicial $k-$chain (with coefficients in $\Z$) is a formal sum of $k$-simplices
\be \sum_{i=1}^N c_i\sigma^i,\mbox{ where }c_i\in \Z,\sigma^i\in{\cal K}\mbox{ is the i-th k-simplex.}\ee

The group of simplicial $k-$chains on $K$ is the free abelian group generated by the sets of all $k$-simplices in $K$, and is denoted by $C_k^\Delta(K)$.
\end{defn}

\begin{defn}[The boundary operator] The boundary operator
\be\partial_k:C^\Delta_k({\cal K})\to C^\Delta_{k-1}({\cal K})\ee
is the homomorphism defined by
\be \partial_k(\sigma)=\sum_{i=0}^k(-1)^i<v^0,\cdots,\hat v^i,\cdots,v^k>,\ee
where the simplex $<v^0,\cdots,\hat v^i,\cdots,v^k>$ is the $i$-th face of $\sigma$ obtained by deleting its $i$-th vertex.
\end{defn}

It is easy to see that
\be \partial_{k}\circ\partial_{k+1}=0\mbox{ for all }k.\ee

\begin{defn}[Cycles and boundaries, homology groups] By (2.9), $Im(\partial_{k+1})\subset Ker(\partial_k)$. The elements of $Im(\partial_{k+1})$ are called boundaries; the elements of $Ker(\partial_k)$ are called cycles. All boundaries are cycles.

The $k$-th simplicial homology group $H_k^\Delta(K,\Z)$ (of coefficients in $\Z$) of a simplicial complex $K$ is the quotient
\be H^\Delta_k(K,\Z)=Ker(\partial_k)/Im(\partial_{k+1}).\ee
\end{defn}

\smallskip
 
 The definition of simplicial homology groups asks some regularity of the space, because we can only define it on a space which admits a triangulation. However, its advantage is that it is relatively easy to understand and calculate, compared to singular homology groups, which we are going to talk about immediately.
 
 \subsubsection{Singular homology (with coefficients in $\Z$)}
 
 \begin{defn}[Standard simplex] The standard $n-$simplex $\Delta_n$ is the convex hull of the points $e_0,\cdots, e_n$ in $\R^n$, where $e_0=(0,\cdots, 0)$ and $e_i=(0,\cdots, 0,1, 0,\cdots, 0)$, the 1 being placed at the $i$-th coordinate.
 \end{defn}
 
 \begin{defn}[Singular simplex] A singular $n-$simplex of a topological space $X$ is a continuous map from $\Delta_n$ to $X$. We denote by $S_n(X)$ the set of all singular $n-$simplices of $X$.
 \end{defn}
 
 \begin{defn}[Singular chain] Let $X$ be a topological space. A singular $k-$chain (of coefficient in $\Z$) on $X$ is a formal sum of $k-$simplices
 \be \sum_{i=1}^Nc_i\sigma^i,\mbox{ where }c_i\in \Z,\sigma^i\in S_k(X)\mbox{ is the i-the singular k-simplex.}\ee
 
 The group of singular $k-$chains on $X$ is the free abelian group generated by $S_k(X)$, and is denoted by $C_k(X,\Z)$.
 \end{defn}
 
 \begin{defn}[The boundary operator] Let $\sigma$ be a singular $n-$simplex of $X$ ($n>0$). The $i$-th face $\sigma_i$ of $\sigma$ is the restriction of the map to the standard $n-$simplex, the convex hull of points $e_0,\cdots, e_{i-1},e_{i+1},\cdots, e_n$.
 
The boundary $\partial\sigma$ of $\sigma$ is defined by $\sum_{i=0}^n(-1)^i\sigma_i$. The boundary of a point (a 0-simplex) is set to be 0. The boundary operator is extended to chains by linearity. So the operator $\partial_n$ sends $C_n(X,\Z)$ to $C_{n-1}(X,\Z)$ (if $n=0, C_{-1}(X,\Z)=0$).
 \end{defn}
 
 Notice that
 \be\partial_n\circ\partial_{n+1}=0\mbox{ for all }k.\ee
 Thus $Im(\partial_{n+1})\subset Ker(\partial_n)$. 
 
 \begin{defn}[Cycles and boundaries, homology groups] Elements in $Im(\partial_{n+1})$ are called boundaries, and elements of $Ker(\partial_n)$ are called cycles. The $k$-th singular homology group $H_k(X,\Z)$ (with coefficients in $\Z$) of a topological space $X$ is the quotient
\be H_k(X,\Z)=Ker(\partial_k)/Im(\partial_{k+1}).\ee
 \end{defn}
 
 Thus we associate to all topological space a sequence of abelian groups.
 
\subsubsection{The relation between these two homology groups}

The singular homology groups are well defined for all topological spaces. But in general it is difficult to calculate directly by definition. So it will be helpful to know, on which kind of spaces these two kinds of homology groups coincide, so that we can decide the singular homology through calculating the simplicial homology, which is much easier.

\smallskip 

First of all, notice that if $K$ is a simplicial complex, there exists a canonical homomorphism from $H_n^\Delta(K,\Z)$ to $H_n(|K|,\Z)$, induced by the map between chains $C_n^\Delta(K,\Z)\to C_n(|K|,\Z)$, which maps each $n$-simplex of $K$ to its characteristic map $\sigma:\Delta^n\to|K|$.

\begin{thm}[c.f.\cite{Hat}, Thm 2.27] The above homomorphism $H_n^\Delta(K,\Z)\to H_n(|K|,\Z)$ is an isomorphism for all $n$.
\end{thm}

After this theorem, If $X$ is a topological space which admits a triangulation $|K|\cong X$, then for any other triangulation $|K'|\cong X$ we have
\be H_n^\Delta({ K},\Z)\cong H_n(|{ K}|,\Z)\cong H_n(X,\Z)\cong H_n(|{ K}'|,\Z)\cong H_n^\Delta({ K}',\Z),\ee
thus we can define simplicial homology groups on X
\be  H_n^\Delta(X,\Z)=H_n^\Delta({K},\Z).\ee

Therefore for every space that admits a triangulation, we can calculate the singular homology by the simplicial homology on this triangulation. The next step is to decide which spaces admit triangulations.

\begin{defn}[$\mu$-smooth triangulation] Let $\mu\ge 1$ be an integer, or $\mu=\infty$. Let $M$ be a $C^\mu$ manifold of dimension $n$. A triangulation $(K,\pi)$ of $M$ is said to be $\mu$-smooth if for every $n-$simplex $\sigma$ of $K$, there exists a chart $(U,\xi)$ of $M$  ($\xi:U\cong V\subset \R^n$) such that $\pi(\sigma)\subset U$, and $\xi\circ \pi$ is affine on $\sigma$. If moreover $K$ is a simplicial $n$-complex, the triangulation is called a $\mu-$smooth $n-$triangulation.
\end{defn}

\begin{thm}[c.f.\cite{Whi},Chap. IV,$\S$ 14B, Thm 12] Every $C^\mu$ manifold $M$ admits a $\mu$-smooth triangulation.
\end{thm}

\begin{rem}Notice that a triangulation $(K,\pi)$ of a $n-$dimensional manifold is automatically a $n$-triangulation.
\end{rem}

\begin{defn}[Smooth simplicial $k-$chain] Let $k\le n$ be two integers, $M$ a $n-$dimensional smooth manifold. For each $\infty-$smooth (we will call it smooth for short) triangulation $(K,\pi)$ of $M$, the image $\Gamma=\pi(\sigma)=\sum_{i=1}^N c_i\pi(\sigma^i)$ of a simplicial $k-$chain $S=\sum_{i=1}^N c_i\sigma^i$ on $K$ under $\pi$  is called a smooth simplicial $k-$chain on $M$. For $d\le k$, the image $\pi(\sigma)$ of a $d-$face $\sigma$ of $K$ is called a $d-$face of $\Gamma$. The boundary $\partial \Gamma$ of $\Gamma$ is defined by
\be\partial\Gamma=\sum_{i=1}^N c_i\pi(\partial\sigma^i)=:\pi(\partial S).\ee
This is a smooth $k-1$-chain.
\end{defn}

\subsection{Transversality}
Intuitively, the property of transversality describes a general position of the intersection of subspaces or submanifolds. It is somehow the opposite of the notion ``tangence".

\begin{defn}Let $M, N$ be two smooth manifolds of dimension $m$ and $n$, $f:M\to N$ a smooth map. Let $\Gamma\subset N$ be a smooth submanifold of $N$. We say that $f$ is transversal to $\Gamma$, denoted $f\pitchfork\Gamma$, if for all $x\in M$ such that $y=f(x)\in\Gamma$, the tangent plane $T_y\Gamma$ and the image $f_*(T_xM)$ generate $T_yN$.
\end{defn}

A most useful property of the transversality is the following proposition.

\begin{pro}[c.f.\cite{Bre} , Chapt II, Thm15.2] Let $M,N,f,\Gamma$ be as in Definition 2.28. If $f\pitchfork\Gamma$, then $f^{-1}(\Gamma)$ is a regular submanifold (also called a embedded submanifold). Moreover the codimension of $f^{-1}(\Gamma)$ in $M$ is the same as that of $\Gamma$ in $N$.
\end{pro}

\begin{cor}Let $M,N,f,\Gamma$ as in the previous proposition. If moreover $M,N,\Gamma$ are orientable, then $f^{-1}(\Gamma)$ is also orientable.
\end{cor}

\nd Denote by $\Gamma'=f^{-1}(\Gamma)$, and $d$ the codimension of $\Gamma$ in $N$, and equivalently the codimension of $\Gamma'$ in $M$, by Proposition 2.29.

$\Gamma$ is orientable, that is, there exists a smooth normal $d-$vector field $v:\Gamma\to\wedge_d(TN)$ on $\Gamma$ which does not vanish, where $v(x)\in \wedge_d(T_x\Gamma^\perp)$. 

Now for each $y\in\Gamma'$, $x=f(y)\in\Gamma$. Since $f$ is transverse, $f_*(T_yM)\oplus T_x\Gamma$ generates $T_xN$. Denote by $N_y\Gamma'$ the orthogonal space of $T_y\Gamma'$ in $T_yM$, and $N_x\Gamma$ the orthogonal space of $T_x\Gamma$ in $T_xN$. We know that $f_*(T_y\Gamma')\subset T_x\Gamma$, hence $f_*(N_y\Gamma')\oplus T_x\Gamma$ generates $N_x\Gamma\oplus T_x\Gamma$. Denote by $\pi_x$ the orthogonal projection of $T_xN$ on $N_x\Gamma$, then $\pi_x\circ f_*(N_y\Gamma')$ is linear and surjective on $N_x\Gamma$. By Proposition 2.29, $N_y\Gamma'$ and $N_x\Gamma$ are both of dimension $d$, hence $\pi_x\circ f_*$ is bijective. Therefore for each $y\in \Gamma'$, set $v'(y)\in N_y\Gamma'$ to be the pre-image of $v(x)$ by $\pi_x\circ f_*$, then $v'$ is a smooth normal $d-$vector field, everywhere non-zero on $\Gamma'$. This gives that $\Gamma'$ is orientable.\qed

\begin{rem}By the proof of Corollary 2.30, the orientation $v'$ of $f^{-1}(\Gamma)$ is called the orientation induced by $f$ from the orientation $v$ of $\Gamma$, denoted by $f^*(v)$. It is easy to see that if $\Gamma$ is a smooth orientable manifold with boundary, $f$ transversal to $\Gamma$ and $\partial\Gamma$, then $\partial f^{-1}(\Gamma)=f^{-1}(\partial\Gamma)$, and moreover if we denote by $\partial v$ the orientation of $\partial \Gamma$ induced by $v$, $\partial f^*(v)$ the orientation of $\partial f^{-1}(\Gamma)$ induced by $f^{-1}(\Gamma)$, then
\be \partial f^*(v)=f^*(\partial v).\ee
\end{rem}

The following theorem says that the subset of all functions transversal to a given submanifold is very large in the set of all smooth functions.

\begin{thm}[Theorem of transversality, c.f.\cite{Hir}, Chapt 3, Thm 2.1]Let $M$ and $W$ be two smooth manifolds. Let $F$ be a closed subset of $M$. Let $Z$ be a smooth submanifold of $W$. Let $f$ be a smooth map from $M$ to $W$ which is transversal to $Z$ at every point of $F$. Denote by $X$ the set of all smooth maps from $M$ to $W$ that coincide with $f$ on $F$. Then the subset of all smooth maps from $M$ to $W$  which coincide with $f$ on $F$ and which are transversal to $Z$ is the intersection of countably many open dense subsets of $X$ with respect to the Whitney topology.
\end{thm}

\begin{rem}Here we do not need to give the precise definition of Whitney topology. All that we have to know is that this is a topology on the space of $C^\infty$ maps, and for all compact set $K\subset M$, for all $\e>0$, for all smooth map $f:M\to W$, the subset of smooth maps $g:M\to W$ such that $|g-f|<\e$ on $K$ is open under Whitney's topology.
\end{rem}


\subsection{The pre-image of a smooth chain under a transverse map}

We want to generalize Proposition 2.29 to all smooth simplicial $k-$chain. Notice that for each $d-$face $\sigma$ of $\Gamma$, the interior $\sigma^\circ$ of $\sigma $ is a smooth regular $d$-submanifold.

\begin{defn}Let $m,n$ be two integers, $M,N$ two smooth manifolds of dimension $m$ and $n$. Let $\Gamma\subset N$ be a smooth $k-$chain ($k<\min\{m,n\}$). We say that a smooth map $f:M\to N$ is transversal to $\Gamma$, if for all $d\le k$, for each $d-$face $\sigma$ of $\Gamma$, $f$ is transversal to the submanifold $\sigma^\circ$. 
\end{defn}

\begin{pro}
 Let $M,N$ be two smooth oriented manifolds of dimension $m$ and $n$, $f:M\to N$ is smooth and proper, $k\le\min\{m,n\}$. Let $\Gamma$ be an smooth simplicial $n-k-$chain in $N$, and suppose that $f$ is transversal to $\Gamma$. Then there exists a triangulation on $M$, under which $f^{-1}(\Gamma)$ is also a smooth simplicial $m-k-$chain in $M$. Moreover,
\be \partial f^{-1}(\Gamma)=f^{-1}(\partial\Gamma).\ee
\end{pro}

\nd 

For each simplicial chain $\gamma$, denote by $|\gamma|$ its support.

Denote by $\Gamma_d$ the union of the interiors of all $d-$faces of $\Gamma$ for $d>0$, and $\Gamma_0$ the union of all $0-$faces. Every $\Gamma_d$ is a $d-$dimensional submanifold. The support $|\Gamma|$ of $\Gamma$ is the disjoint union of all these $\Gamma_d$. By hypothesis, $f$ is transversal to all these $\Gamma_d$, hence by Proposition 2.29, $f^{-1}(\Gamma_d)$ is a smooth submanifold of dimension $m-n+d$, and the support of $f^{-1}(\Gamma)$ is their disjoint union.

We do a triangulation $K'$ on $M$, such that the restriction to each manifold $f^{-1}(\Gamma_d)$ is a finite triangulation on it. This kind of triangulation exists because we just have to find a triangulation on $M$ which is transversal to each $f^{-1}(\Gamma_d)$. The finiteness of the triangulation on each $f^{-1}(\Gamma_d)$ comes from the fact that $f$ is proper. Thus the intersection of $K'$ with each $f^{-1}(\Gamma_d)$ is a decomposition of it into smooth polygons. We do subdivisions on these polygons so that it becomes a triangulation. Thus we get a triangulation on each $f^{-1}(\Gamma_d)$. Last we subdivide every face of $K'$ that meets $f^{-1}(\Gamma_d)$, and get a finer triangulation $K$ on $M$, under which each $f^{-1}(|\sigma|)$ is the support of a simplicial chain, for each face $\sigma$ of $\Gamma$. Denote also by $f^{-1}(\sigma)$ this simplicial chain, that is, the sum of all the simplices contained in it, with the proper orientation.

Now we work with this triangulation. By linearity, all we have to prove is that (2.37) is true if $\Gamma$ is a simplicial simplex.

We know that $|f^{-1}(\Gamma)|$ is a stratified space, with layers $f^{-1}(\Gamma_d)$. So to prove (2.37), we only have to prove that
for each $d$, for each $x\in\Gamma_d$, and each $y\in f^{-1}(x)$, the local structure of $f^{-1}(\Gamma)$ around $y$ is the same as that of $\Gamma$ around $x$.  Without loss of generality, we can even suppose that $\Gamma$ is a simplex in $\R^n$, and $M$ is a ball in $\R^m$, since we only want to get local structure. We suppose also that $x$ and $y$ are the origins.

So denote by $\sigma$ the $d-$face of $\Gamma$ such that $0\in\sigma^\circ$, and for each $l>0$ there are $k_l$ $d+k-$faces of $\Gamma$ who contain $\sigma$. Denote by $P$ the $d-$plane containing $\sigma$, and $P^\perp$ the orthogonal space of $P$ in $\R^n$. Take a small enough ``cylinder'' $B=B_d\times B_{n-d}\subset\R^n$, where $B_d$ and $B_{n-d}$ are balls contained in $P$ and $P^\perp$ respectively with center $0$, such that $B$ only meets faces of $\Gamma$ that contain $\sigma$. Thus $B\cap\Gamma=B_d\times C$, where $C$ is a cone in $B_{n-d}$, with $k_l$ $k-$faces.

Now since $f$ is transversal to $\Gamma$, if $d<n-m$, then the image of $f$ will not meet $\sigma$, and there is nothing to prove. So suppose $d\ge n-m$. Denote by $\pi^\perp$ the orthogonal projection of $\R^n$ to $P^\perp$, and set $g=\pi^\perp\circ f:M\to \sigma^\circ$. By transversality, $df(0)(\R^m)\oplus P=\R^n$. Hence $df(0)(\R^m)$ contains $P^\perp$.  Therefore $dg(0):\R^m\to\R^{n-d}$ is onto. Thus by Theorem 7.3 of Chapt II of \cite{Bre}, there are diffeomorphisms $\phi$ of a neighborhood $U\subset M$ of $y=0$ and $\psi$ of a neighborhood $V\subset \sigma^\circ$ of $x=0$ such that the function $\phi^{-1}\circ g\circ\psi$ is the projection of $\R^m$ onto $\R^{n-d}:(x_1,\cdots,x_m)\mapsto(x_1,\cdots,x_{n-d})$.

As a result, locally, $g^{-1}(C)$ is diffeomorphic to $C\times \R^{m-n+d}$. But $\Gamma$ is locally $C\times\R^d$, hence $f^{-1}(\Gamma)$ is locally diffeomorphic to $C\times \R^{m-n+d}$. Note that under this homeomorphism, $f^{-1}(\sigma)$ is $\{0\}\times\R^{m-n+d}$ around $y$, thus the structure around $y$ of $f^{-1}(\Gamma)$ is the same as that around $x$ in $N$. \qed

\section{Topological minimal set}

We begin now to generalize the notion of MS-minimizers to higher codimensions.

In this section, let $U\subset\R^n$ be an open set, and for every oriented compact $d-$manifold with boundary $\Gamma\subset U$, $\Gamma$ represents also the chain $\sum_{i=1}^N\Sigma^i$, where $\{\Sigma^i,1\le i\le N\}$ is the set of all smooth simplicial $k-$simplices of no matter which smooth triangulation on the manifold with boundary $\Gamma$ which preserve the orientation. This is well defined, since by Theorem 2.20, it does not depend on triangulations.

\begin{defn}[Topological competitors] Let $E$ be a closed set in $\R^n$. We say that a closed set $F$ is a topological competitor of dimension $d$ ($d<n$) of $E$, if there exists a ball $B\subset\R^n$ such that

1) $F\bs B=E\bs B$;

2) For each Euclidean $n-d-1$-sphere $S\subset\R^n\bs(B\cup E)$, if $S$ represents a non-zero element in the singular homology group $H_{n-d-1}(\R^n\bs E;\Z)$, then it is also non-zero in $H_{n-d-1}(\R^n\bs F;\Z)$.

We also say that such a $F$ is a topological competitor of $E$ with respect to the ball $B$.
\end{defn}

Note that $\R^n\bs E$ and $\R^n\bs F$ are open in $\R^n$, hence they are smooth $n-$manifolds, on which the singular homology and simplicial homology coincide. So 
it is equivalent to replace the condition 2) by 

2') For each Euclidean $n-d-1$-sphere $S\subset\R^n\bs(B\cup E)$, if $S$ represents a non-zero element in the simplicial homology group $H^\Delta_{n-d-1}(\R^n\bs E;\Z)$, then it is also non-zero in $H^\Delta_{n-d-1}(\R^n\bs F;\Z)$.

\begin{rem}
 For the case $d=n-1$, a one-dimensional sphere is composed of two points. And one can check that the sphere is zero in homology if and only if the two points are in the same connected component. Hence in the case of codimension 1, the definition of topological competitors coincides with that of MS-competitors.
\end{rem}

\begin{rem}
 In Definition 3.1, if $F$ is a topological competitor of $E$, and $B$ is the associated ball which verifies 1) and 2), then for all ball $B'\supset B$, it verifies 1) and 2), too; on the contrary for a smaller ball $B'$ that contained in $B$ who verifies 1), the condition 2) is not necessarily true for $B'$. See the picture 3-1 for a counterexample, where the set $F$ is a competitor of $E$ with respect to the larger ball $B$, but not for the smaller one $B'$. In fact, $E$ separates the two points (the two little crosses in the picture), but $F$ does not. 
 
 \centerline{\includegraphics[width=0.45\textwidth]{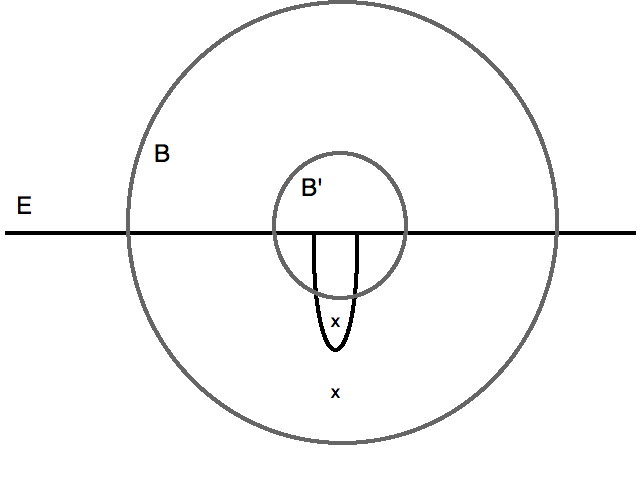} \includegraphics[width=0.45\textwidth]{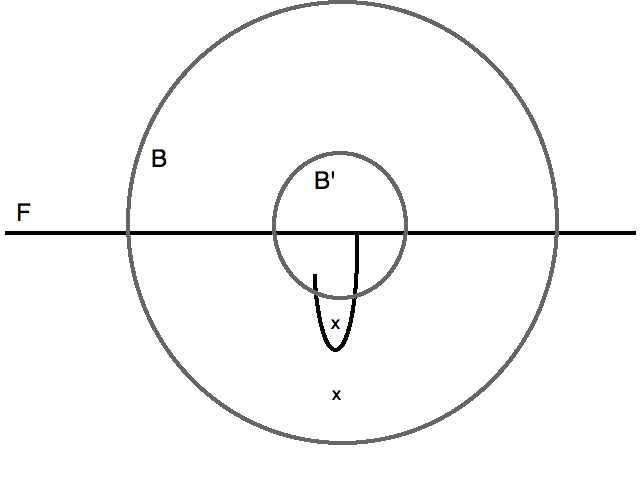}}
\nopagebreak[4]
\centerline{3-1}
\end{rem}

\begin{defn}[Topological minimal set] A topological minimal set is a minimal set defined as in Definition 1.1 where we take the topological competitors as the competitor class.
\end{defn}

By Remark 3.2, in codimension 1, the definition of topological minimal sets coincides with that of MS-minimal sets.

We also have the following local version:

\begin{defn}[locally topological minimal sets] We say that a closed set $E$ is $d-$dimensional locally topologically minimal in a ball $B$ if $E$ minimizes the $d-$Hausdorff measure among all its competitors with respect to the ball $B$.
\end{defn}

\begin{rem}
 It is also possible to define the topological competitors by the homotopy groups $\pi_{n-d-1}$, which is somehow more natural, and the definition of $\pi_0$ coincides also with the separation condition when $d=n-1$. But we know a simple example in $\R^3$, a Fox-Artin arc, which admits a deformation that does not preserve the group $\pi_1$ of its complement. See \cite{Artin} for more detail.
\end{rem}

However we are lucky with our definition by homology, as stated in the following proposition.

\begin{pro} Let $E\subset\R^n$ be closed. Then every Almgren competitor $F$ of $E$ is a topological competitor of $E$ of dimension $d$, for any $d<n$.
\end{pro}

\nd Let $E$ be a closed set in $\R^n$, and $F$ be an Al-competitor of $E$, that is, there exists a ball $B=B(x,R)\subset\R^n$, and a Lipschitz map $f$ from $\R^n$ to $\R^n$ such that
\be f(B)\subset B,f|_{\R^n\bs B}=id\mbox{ and }f(E)=F.\ee

Fix a $d<n$. For each $r>0$, set $rB=B(x,rR)$. We will prove that $F$ is a topological competitor of $E$ with respect to $2B$.

Let $S\subset\R^n\bs(2B\cup E)$ be a euclidean $n-d-1-$sphere, then $S\subset\R^n\bs F$. Suppose that there exists a smooth simplicial $n-d-$chain $\Gamma\subset\R^n\bs F$ whose boundary is $S$, we want to show that there exists also a smooth simplicial $n-d-$chain $\Gamma'\subset\R^n\bs E$ with $S$ as its boundary.

Set $\d=\mbox{dist}(\Gamma,F)$. Then we can find a smooth map $g$ from $\R^n$ to $\R^n$, with
\be g_{\R^n\bs\frac32 B}=f=id,\ee
and
\be ||g-f||_\infty=\sup\{||f(x)-g(x)||,s\in\frac32 B\}<\frac\d4.\ee

By (3.9), $g$ is transversal to $\Gamma$ on $\R^n\bs\frac32 B$, since it is identity.

By Theorem 2.33, the set of smooth fonctions from $\R^n$ to $\R^n$ which are transversal to $\Gamma$ and which coincide with $g$ on $\R^n\bs\frac32 B$ is a $G_\delta$-dense subset in the set of all smooth functions from $\R^n$ to $\R^n$ which coincide with $g$ on $\R^n\bs\frac 32 B$, under the topology of Whitney (see Remark 2.34). In particular, there exists a smooth map $h$ from $\R^n$ to $\R^n$ such that
\be h_{\R^n\bs\frac32 B}=g=id,\ee
\be ||g-h||_\infty=\sup\{||h(x)-g(x)||,x\in\frac32 B\}<\frac\d4,\ee
and
\be h\mbox{ is transversal to }\Gamma.\ee

Note that (3.11) means that $h$ is surjective and proper. By Proposition 2.36, (3.13) gives that $\Gamma'=h^{-1}(\Gamma)$ is a smooth simplicial $k-$chain in $\R^n$, and moreover
\be \partial\Gamma'=f^{-1}(\partial\Gamma)=S.\ee

But $\Gamma'\cap E=\emptyset$. In fact, if $x\in\Gamma'\cap E$, then $f(x)\in F$, and $h(x)\in\Gamma$. But $|h(x)-f(x)|\le||g-h||_\infty+||g-f||_\infty\le\frac\d 2$, thus
\be d(F,\Gamma)\le\frac\d 2<\d=\mbox{dist}(\Gamma,F),\ee
which leads to a contradiction.

Thus we have obtained a smooth simplicial chain $\Gamma'$ in $\R^n\bs E$ whose boundary is $S$. 

Therefore, $F$ is a topological competitor of $E$ with respect to $2B$.\qed 

\begin{rem}
 The proof can be trivially generalized to prove that, if $S$ is a smooth $n-d-1-$surface without boundary, and is non zero in $H_{n-d-1}(\R^n\bs E)$ for some $d-$dimensional set $E$, and if $f:\R^n\to\R^n$ is a map that is different from the identity only in some compact set $B$ that does not touch $S$, then $S$ is also non zero in $H_{n-d-1}(\R^n\bs f(E))$.
\end{rem}

\begin{cor} All topological minimal sets are Al-minimal sets.
\end{cor}

Next we show that the topological minimality does not depend on the ambient dimension.

\begin{pro} Let $E$ be a $d-$dimensional set. Then if there exists $m>d$ such that $E$ is topological minimal in $\R^m$, $E$ is topological minimal in $\R^n$ for all $n>d$ such that $E\subset\R^n$.
\end{pro}

\nd Let $E\subset\R^n$. We will just show the proposition for two cases: 1) $E$ is topological minimal in $\R^n$ implies that it is also topological minimal in $\R^{n+1}$; 2) $E$ is topological minimal in $\R^{n+1}$ implies that $E$ is topological minimal in $\R^n$. The proposition follows then by induction.

1) Let $E\subset\R^n$ be a topological minimizer of dimension $d$. Let $F$ be a topological competitor of $E$ in $\R^{n+1}=\R^n\times \R$. That is, there exists a $n+1$-dimensional ball $B\subset\R^{n+1}$, such that $F\bs B=E\bs B$, and for each $n-d-$sphere $S\subset\R^{n+1}\bs(B\cup E)$, if $S$ is non-zero in $H_{n-d}(\R^{n+1}\bs E,\Z)$, then it is non-zero in $H_{n-d}(\R^{n+1}\bs F,\Z)$.

By Remark 3.3, we can suppose that the center of $B$ belongs to $\R^n$. Notice that $E\subset\R^n$, hence $F\bs B=E\bs B\subset\R^n$. Denote by $\pi$ the orthogonal projection on $\R^n$, then by Proposition 3.7, $\pi(F)\subset\R^n$ is a topological competitor of $E$ in $\R^{n+1}$, with respect to the ball $2B$. Since $\pi$ is 1-Lipschitz, $\pi(F)$ is a better topological competitor of $E$ than $F$, that is
\be H^d(\pi(F)\bs F)\le H^d(F\bs\pi(F)).\ee

Denote by $B'=\pi(2B)=2B\cap\R^n$ (since the center of $2B$ belongs to $\R^n$) a $n-$dimensional ball contained in $\R^n$. We want to show that $\pi(F)$ is a topological competitor of $E$ in $\R^n$ with respect to $B'$.

Let $S\subset\R^n\bs(B'\cup E)$ be a $n-d-1-$sphere. If $S$ is zero in $\R^n\bs\pi(F)$, there exists a simplicial $n-d$-chain $\Gamma\subset\R^n\bs\pi(F)$, such that $\partial\Gamma=S$.

Denote by $U\subset\R^n\bs B'$ the $n-d$ ball in $\R^n$ such that $\partial U=S$.

We know that $H_{n-d}(\R^n,\Z)=0$ for $d\ge 1$, and that $U-\Gamma$ is a simplicial $n-d$-chain without boundary, therefore there exists a $n-d+1$-chain $R\subset\R^n$ such that $\partial R=U-\Gamma$.

Set $T=R\times\{1\}-R\times\{-1\}+\Gamma\times[-1,1]$ in $\R^{n+1}$, then $T$ is a simplicial $n-d+1$-chain, with $\partial T=S\times[-1,1]+U\times\{1\}-U\times\{-1\}$. Thus $\partial T$ is a compact topological manifold without boundary, and is smooth except near $S\times\{-1,1\}$. So we move $\partial T$ a little in a small neighborhood of $S\times\{-1,1\}$, and obtain that in $\R^{n+1}\bs(2B\cup\pi(F))$, it is homotopic to the $n-d$-sphere $S'$ whose center is in $\R^n$ and is such that $S'\cap\R^n=S$.

But $\partial T$ is zero in $H_{n-d}(\R^{n+1}\bs\pi(F),\Z)$, hence $S'\subset\R^{n+1}\bs 2B$ is zero in $H_{n-d}(\R^{n+1}\bs\pi(F),\Z)$. Recall that $\pi(F)$ is a topological competitor of $E$ in $\R^{n+1}$ with respect to $2B$, hence $S'$ is zero in $H_{n-d}(\R^{n+1}\bs E,\Z)$. That is, there exists a smooth simplicial $n-d+1$-chain $R'\subset\R^{n+1}\bs E$ whose boundary is $S'$. Take a ball $B_1\subset\R^{n+1}$ such that $R'\subset B_1$.

Denote by $i$ the embedding from $\R^n\to\R^{n+1}$. Then $i$ is transversal to $R'$ on $S'\cup B_1^C$. Thus by Theorem 2.33, the set of all smooth functions from $\R^n$ to $\R^{n+1}$ which coincide with $i$ on $S'\cup B_1^C$ and which are transversal to $R'$ is $G_\d$-dense in the set of all smooth functions in $C^\infty(\R^n,\R^{n+1})$ which coincide with $i$ on $S'\cup B_1^C$. In particular, since the set of embeddings is open (\cite{Hir} Chapt 2, Thm 1.4), there exists an embedding $g\in C^\infty(\R^n,\R^{n+1})$, transversal to $R'$, such that $||g-i||_\infty$ is small enough such that $g^{-1}(R')\cap E=\emptyset$. As a result, Proposition 2.36 gives that $\Gamma'=g^{-1}(R')$ is a smooth simplicial $n-d-$chain in $\R^n$ with $\partial\Gamma'=g^{-1}(S')=S$. Hence $S$ is zero in $H_{n-d-1}(\R^n\bs E,\Z)$.

Therefore $\pi(F)$ is a topological competitor of $E$ in $\R^n$. By the topological minimality of $E$, we have
\be H^d(E\bs\pi(F))\le H^d(\pi(F)\bs E).\ee

Combine with (3.19), we get
\be H^d(E\bs F)\le H^d(F\bs E),\ee 
and hence $E$ is topological minimal in $\R^{n+1}$.

2) Let $E\subset\R^n$ be a $d-$dimensional topological minimal set in $\R^{n+1}=\R^n\times\R$. Let $F\subset\R^n$ be a topological competitor of $E$ in $\R^n$. That is, there exists a ball $B\subset\R^n$ such that $F\bs B=E\bs B$, and for all $n-d-1-$sphere $S\subset\R^n\bs(B\cup E)$, if $S$ is the boundary of a smooth simplicial $n-d-$chain in $\R^n\bs E$, then it is also the boundary of some smooth simplicial $n-d$-chain in $\R^n\bs F$.

Denote by $B'$ the ball in $\R^{n+1}$ with the same center and radius as $B$. Then $F\bs B'=E\bs B'$. We want to show that $F$ is also a topological competitor of $E$ in $\R^{n+1}$.

Let $S'\subset\R^{n+1}\bs(F\cup B')$ be a $n-d$-dimensional sphere, which is zero in $H_{n-d}(\R^{n+1}\bs F,\Z)$. We would like to show that $S'$ is zero in $H_{n-d}(\R^{n+1}\bs E,\Z)$.

If $S'\cap\R^n=\emptyset$, then since $S'$ is connected, and $\R^n$ separates $\R^{n+1}$ into two components, $S'$ has to be contained in one of these two components, while both are homeomorphic to $\R^{n+1}$, so that their $n-d$ homology group $H_{n-d}$ is trivial. As a result, $S'$ is zero in $H_{n-d}(\R^{n+1}\bs E,\Z)$;

If $S'\cap\R^n\ne\emptyset$, but the intersection is just a point, then this is the case when $\R^n$ is tangent to $S'$. The argument is similar to the above one;

If $S'\cap\R^n\ne\emptyset$ and the intersection is non trivial, then the intersection $S=S'\cap \R^n$ is a $n-d-1$-sphere. Since $S'$ is zero in $H_{n-d}(\R^{n+1}\bs F,\Z)$, there exists a smooth simplicial $n-d+1$-chain $\Gamma\subset\R^{n+1}\bs F$ whose boundary is $S'$. Then by an argument similar to the argument in 1) which has been used to show that $S$ is a zero element, there exists a smooth simplicial $n-d$-chain $R$ in $\R^n\bs E$ whose boundary is $S$. Take $D\subset\R^n\bs B'$ a compact smooth manifold with boundary of dimension $n-d$ whose boundary is $-S$. Then $D+R$ is a $n-d-$chain of boundary zero in $\R^n$. Then there exists a $n-d+1-$chain $U\subset\R^n$ such that $\partial U=D+R$, since $H_{n-d}(\R^n,\Z)$ is trivial.

Set $R'=R\times[-1,1]+U\times\{1\}-U\times\{-1\}$, then $R'\cap E=\emptyset$ and $\partial R'=S\times[-1,1]+D\times\{1\}-D\times\{-1\}$, which is homotopic to $S'$. Hence $S'$ is zero in $H_{n-d}(\R^n\bs E,\Z)$.

As a result, $F$ is a topological competitor of $E$ in $\R^{n+1}$ with respect to $B'$, which gives
\be H^d(E\bs F)\le H^d(F\bs E).\ee

Hence $E$ is topologically minimal in $\R^n$.\qed

Next we are going to discuss the product of two sets. That is, if $E=E_1\times E_2$, with $E_i\subset\R^{n_i}$, $E\subset\R^{n_1+n_2}$, then what is the relation among minimalities of the three?

One direction is simple, that is , if $E$ is Al-minimal in $\R^{n_1+n_2}$, then the $E_i$'s are also Al-minimal in $\R^{n_i},i=1,2$. (c.f.\cite{XY10}, section 12).

For the other direction, if the two $E_i$'s are Al-minimal in $\R^{n_1},i=1,2$, we do not know yet how to prove the Al-minimality of their product. In fact, even for the case when one of the $E_i$ is $\R$, we do not know how to prove it. A very natural way to deal with products is always the ``slicing'' method. Take $E_1\times\R$ for example, we would like to use the minimality of $E_1$, hence we look at the slice $\pi^{-1}(x)$ for each $x\in\R$, where $\pi$ denotes the orthogonal projection on $\R$. If we could show that for any deformation $F$ of $E_1\times\R$, each slice $\pi^{-1}(x)\cap F$ is a deformation of $E_1$, then we would win because $E_1$ is minimal. But unfortunately this is not true. 

For the general case when $E=E_1\times E_2$ the situation is more complicated, because in general the slices $\pi^{-1}(x)\cap F,x\in E_2$ for a deformation $F$ of $E_1\times E_2$ are not even of the same dimension as $E_1$, because $E_2$ is not of full dimension.

Hence probably a more intelligent way of slicing is needed.

The above is the case for Al-minimal sets. Here for our topological minimal sets, at least we know how to prove the topological minimality of the product of a topological minimal set with $\R^n$.

\begin{pro}
 Let $E\subset\R^n$ be a topological minimal set of dimension $d$ in $\R^n$, then for all $m>0$, $E\times\R^m$ is topologically minimal of dimension $d+m$ in $\R^{n+m}$.
\end{pro}
 
\nd We will only prove it for $m=1$. The case for $m>1$ comes then by induction.

First we restrict the class of topological competitors of $E\times\R$ to the subset of all rectifiable topological competitors. In fact, if $F$ is a non-rectifiable topological competitor of $E\times\R$, with respect to a ball $B$ for example, then the part inside of $B$ is not rectifiable, because outside $B$, $F=E\times\R$, and $E$ is rectifiable since it is topologically minimal and hence Al-minimal by Corollary 3.17. (For the rectifiability of Al-minimal sets, see for example \cite{DS00} 2.11.)

But $F$ is not rectifiable in $B$ means that $F$ is not Al-minimal in $B$ (c.f.\cite{DS00} Thm 2.11). Thus there exists an Al-competitor $F'$ of $F$ in $B$ such that $0<\d=H^d(F\cap B)-H^d(F'\cap B)$. By Proposition 3.7, $F'$ is also a topological competitor of $E\times\R$ in $B$. Denote by $rB$ the ball with the same center of $B$ and whose radius is $r$ times of that of $B$.  Then by Lemma 5.2.2 of \cite{Fv}, there exists an $n+1$ dimensional complex $S$ (whose support is denoted by $|S|$), with $\overline{\frac32 B}\subset|S|^\circ$ and $|S|\subset 2B$, and a deformation $f$ in $2B$, such that $f|_{|S|}$ is some Federer-Fleming projection on $S$, $f(\overline B)\subset \frac43 B$, and $H^d(f(F')\cap 2B)< H^d(F'\cap 2B)+\frac\d2< H^d(F\cap 2B)$. Then still by Proposition 3.7, $F''=f(F')$ is a topological competitor of $E$ in $2B$. On the other hand, we have that $F''\cap \frac 43 B$ is composed of faces of $|S|$, hence is rectifiable, and $F''\bs\frac 43 B$ is contained in $f(F'\bs \overline B)=f(E\times\R\bs\overline B)$, which is rectifiable, hence $F''$ is rectifiable. Thus we get a better rectifiable competitor $F''$ than $F$.

So we have that
\be\begin{split}\mbox{ a set }&F\mbox{ is a topological minimizer of dimension }d\Leftrightarrow F\mbox{ minimizes the }\\
&d-\mbox{dimensional Hausdorff measure among all its rectifiable topological competitors}.
   \end{split}
\ee

So let $F$ be a rectifiable topological competitor of $E\times\R$. By definition, we can find $R>0$ such that 
\be F\bs B(0,R)=(E\times\R)\bs B(0,R),\ee
and that for all $n-d-1$-sphere $S\subset\R^{n+1}\bs((E\times\R)\cup B(0,R)),$
\be \begin{split}S\mbox{ is non-zero in }&H_{n-d-1}(\R^{n+1}\bs E\times\R,\Z)\\
     &\Rightarrow\mbox{it is non-zero in }H_{n-d-1}(\R^{n+1}\bs F,\Z).
    \end{split}\ee

We replace $B(0,R)$ by $C(R):=B_n(0,R)\times[-R,R]$, where $B_n(0,R)=B(0,R)\cap\R^n$. Then (3.25) and (3.26) are still true.

Denote by $\pi$ the orthogonal projection on $\R$, and for all $t\in[-R,R]$, by $F_t$ the slice $\pi^{-1}(t)\cap F$ of $F$. We claim that $F_t-t$ is a topological competitor of $E$ in $\R^n$. 

In fact, let $S\subset\R^n\bs[(F_t-t)\cup B_n(0,R)]$ be a $n-d-1$ sphere, which represents a zero element in $H_{n-d-1}(\R^n\bs(F_t-t),\Z)$. This means that there exists a singular $n-d$-chain $\Gamma\subset\R^n\bs(F_t-t)$ such that $\partial \Gamma=S$. Then $\Gamma_t:=\Gamma\times\{t\}\subset\R^{n+1}\bs F$ is a $n-d$-chain and $\partial\Gamma_t=S\times\{t\}$. Notice that $S\times\{t\}\subset(\R^n\bs[(F_t-t)\cup B_n(0,R)])\times\{t\}=[\R^n\bs B_n(0,R)]\times\{t\}\bs F_t\subset\R^{n+1}\bs(F\cup C(R)),$ and that $\Gamma_t\subset\R^{n+1}\bs F$ is such that $\partial \Gamma_t=S\times\{t\}$, hence $[S\times\{t\}]$ is zero in $H_{n-d-1}(\R^{n+1}\bs F,\Z)$. As a result $[S\times\{t\}]$ is zero in $H_{n-d-1}(\R^{n+1}\bs (E\times\R),\Z)$, since $F$ is a topological competitor of $E\times\R$ with respect to $B(0,R)$. Thus there exists a singular $n-d$-chain $\Gamma'\subset\R^{n+1}\bs(E\times\R)$ such that $\partial\Gamma'=S\times\{t\}$. Denote by $\Gamma''=p(\Gamma')$, where $p$ denotes the projection to $\R^n$. Then $\Gamma''\subset\R^n\bs E$ and $\partial\Gamma''=S$. Hence $S$ is zero in $H_{n-d-1}(\R^n\bs E,\Z)$, and hence $F_t-t$ is a topological competitor of $E$.

By the topological minimality of $E$,
\be\begin{split}H^d(\pi^{-1}(t)\cap F\cap C(R))&=H^d([(F_t-t)\cap B_n(0,R)]\times\{t\})\\
    &=H^d((F_t-t)\cap B_n(0,R))\ge H^d(E\cap B_n(0,R)).
   \end{split}\ee
 
However by the coarea formula (c.f.\cite{Fe} Thm 3.2.22)
\be\int_{F\cap C(R)}||apD\pi(x)||dH^{d+1}(x)=\int_{-R}^RH^d(\pi^{-1}(t)\cap F\cap B_n(0,R)\times\{t\}) dH^1(t).\ee
But $\pi$ is 1-Lipschitz, therefore $||apD\pi(x)||\le 1$, and
\be\int_{F\cap C(R)}||apD\pi(x)||dH^{d+1}(x)\le H^{d+1}(F\cap C(R)).\ee
On the other hand, by (3.27)
\be\int_{-R}^RH^d(\pi^{-1}(t)\cap F\cap B_n(0,R))dH^1(t)\ge 2RH^d(E\cap B_n(0,R)).\ee
Thus we get
\be H^{d+1}(F\cap C(R))\ge 2RH^d(E\cap B_n(0,R))=H^{d+1}(E\cap C(R)),\ee
which gives the topological minimality of $E\times \R$.\qed

\section{Existence theorems}

After the general discussion above, in this section we are going to prove an existence theorem. But instead of considering only minimal sets in the whole $\R^n$, we are going to give the existence result for a little more general setting. 

\begin{defn}
 Let $\Omega\subset\R^n$ be a closed set. For each $\e>0$, denote by $\Omega_\e=\{x:\mbox{dist}(x,\Omega)<\e\}$ the $\e-$neighborhood of $\Omega$. We say that $\Omega$ admits a Lipschitz neighborhood retraction, if there exists $\e>0$, and a Lipschitz map $f$ from $\Omega_\e$ to $\Omega$, such that $f|_{\Omega}=id$.
\end{defn}

\begin{thm} Let $\Omega\subset\R^n$ be a compact set that admits a Lipschitz neighborhood retraction from $\Omega_\e$ to $\Omega$ with $\e>0$. Let $h$ be a continuous function from $\Omega$ to $\R$, and suppose that $1\le h\le M$ on $\Omega$. Let $\{w_j\}_{j\in J}$ be a family of smooth $n-d-1-$ surfaces in $\R^n\bs\Omega_\e$ which are non-zero in $H_{n-d-1}(\R^n\bs \Omega)$. Set
 \be \mathfrak F=\{F\subset\Omega\mbox{ closed, and for all }j\in J, w_j\mbox{ is non-zero in }H_{n-d-1}(\R^n\bs F)\}.\ee
Define
\be J_h(F)=\int_Fh(x)dH^d(x).\ee
Then there exists $F_0\in\mathfrak F$ such that
\be J_h(F_0)=\inf\{J_h(F);F\in\mathfrak F\}.\ee
\end{thm}

\nd The idea is not complicated. Set $m=\inf\{J_h(F);F\in\mathfrak F\}$. Suppose also that $m<+\infty$ because otherwise we have nothing to prove. Given a minimizing sequence $\{F_k\}_{k\in\N}$, that is
\be\lim_{k\to\infty} J_h(F_k)=m,\ee
we can try to extract a subsequence that converges to a set $F_0$ for the Hausdorff distance. Denote still by $F_k$ this converging subsequence. If we can prove that

1) $J_h(F_0)=m;$

2) $F_0\in\mathfrak F$,

then $F_0$ will be a solution, and the theorem will be proved.

For 1), we need the lower semi continuity, that is 
\be J_h(\lim_kF_k)\le\liminf_{k\to\infty}J_h(F_k).\ee

%
%

In fact, the lower semi-continuity (4.7) is clearly not true for all converging sequence. But we do not need that, either. All we have to do is to find a minimizing sequence for which the lower semi-continuity holds.

Now suppose we are given any minimizing sequence $F_k$, that is, (4.6) holds. We would like to modify each $F_k$ a little to a new set $F_k'$, so that after the modification, $F_k'$ is a minimizing sequence, which preserve the topological condition for $w_j$, and in addition the lower semi-continuity holds for this $\{F_k'\}$. But since $\Omega$ is closed, and the modification always has to take place in a neighborhood of each $F_k$, so we need the retract property of the domain $\Omega$.

We know that $\Omega$ is the Lipschitz retract of its $\e-$neighborhood $\Omega_\e$ for some $\e>0$. Denote by $L$ the Lipschitz constant of this retraction.

Next we extend $h$ to $\Omega_\e$ by 
\be h'(x)=\left\{\begin{array}{cc}h(x),&x\in\Omega;\\
                 10L^dM,&x\in  \Omega_\e\bs\Omega. \end{array}
\right.\ee
Then $h'$ is lower semi continuous. Hence there exists a sequence of continuous functions $g_n$, with $g_n\le g_{n+1}$, $g_n|_{\overline\Omega}=h$, and $g_n\uparrow h'$ pointwisely. (For example let $g_n'(x)=\inf\{h'(y)+n|y-x|,y\in \Omega_\e\}$, then $g_n'$ is continuous and converge monotonously to $h'$; let $g$ be a continuous extension of $h$ on $\Omega_\e$ which is smaller than $h'$, and set $g_n=\sup\{g_n',g\}$.)  Denote by $J_{h'}(F)=\int_{F\cap \Omega_\e} h'dH^d$, and $J_{g_n}(F)=\int_{F\cap \Omega_\e} g_ndH^d$ for any closed set $F$. Notice that $J_{h'}(F)=J_{g_n}(F)=J_h(F)$ for all $F\subset\Omega$.

We cover the set $\overline\Omega_{\frac\e3}\bs\Omega_{\frac\e5}$ with a finite set ${\cal D}=\{Q_n,1\le n\le N\}$ of dyadic cubes of size $\frac{\e}{10000}$. Denote by $ G$ the closed set $\Omega_{\frac\e3}\cup(\cup_{n=1}^NQ_n)$. Then $ G\subset\Omega_{\frac\e2}$. And the boundary of $ G$ is a union of faces of dyadic cubes. By Lemma 12.2 of \cite{DS00}, there exists $r>0$ such that $ G$ is a Lipschitz retract of its neighborhood $ G_r=\{x\in\R^n:\mbox{dist}(x, G)<r\}.$

Denote also by $U$ the $\frac{1}{100}\e$ neighborhood $\Omega_{\frac{1}{100}\e}$ of $\Omega$.

Now by Lemma 5.2.2 of \cite{Fv}, for each $k$, we can find an $n-$dimensional complex (a complex composed of polygons) $S_k\subset U$, each polygon of $S_k$ is of size similar to $\frac{\e}{10000}$, and $F_k'\subset U$ such that 

1) The set $F_k'$ is a deformation of $F_k$ in $U$. And if we denote by $|S_k|$ the support of $S_k$, then $\Omega\subset |S_k|^\circ, |S_k|\subset U$, and $F_k'\subset |S_k|$ is a union of $d-$faces of $S_k$.

2) $J_{g_k}(F_k')\le (1+2^{-k})J_{g_k}(F_k)$; 

3) $F_k'$ minimizes $J_{g_k}$ among all deformations of $F_k$ on unions of $d-$faces of $S_k$..

We can also ask the $S_k$ to be uniformly round, that is, the angle between any two faces of any polygon in any $S_k$ is bounded below by a constant which does not depend on $k$. This uniform roundness gives a constant $L'$, such that for any $k$, and any $d-$dimensional set $F$ contained in $|S_k|$, there exists a Federer-Fleming projection from $F$ to the $d-$skeleton of $S_k$ which is $L'$ Lipschitz. See Lemma 4.3.2 of \cite{Fv} for more detail. 

Now the support of $S_k$ and the dyadic cubes in $\cal D$ is relatively far, so by Theorem 2.3 of \cite{Fv}, there exists a complex $S_k'$, such that $S_k'$ contains $S_k$ and all the cubes in $\cal D$, $|S_k'|= G$, and $S_k'$ has the same roundness as $S_k$. Hence the $S_k'$ are uniformly round, too.

We take $F_k''$ the set which minimize $J_{g_k}$ among all deformations of $F_k$ on the union of $d-$faces of $S_k'$. This minimizer exists. In fact, for any proper closed subset of any $d-$face of $s_k$, we can retract it to the boundary of this $d-$face by a radial projection on the $d-$face. Hence for any subset of the $d-$skeleton of $S_k'$, we can deform it to a subset of itself which is a union of $d-$faces of $S_k'$. Thus the minimizer is actually taken over all the unions of some $d-$faces of $S_k'$ that is a deformation of $F_k$. Then since the set of unions of $d-$faces of $S_k'$ is finite, a minimizer exists, and is a union of some $d-$faces of $S_k'$. Moreover, since $F_k''$ minimizes $J_{g_k}$, we have
\be J_{g_k}(F_k'')\le J_{g_k}(F_k')\le(1+2^{-k})J_{g_k}(F_k).\ee

We want to prove that the $F_k''$ are uniformly quasiminimal in $\Omega_\e$, which will give the uniform concentration property for $F_k''$, and thus implies the lower semi continuity of Hausdorff measure. We give first the definition of quasiminimality.

\begin{defn}[Quasiminimality] Let $U\subset\R^n$ be an open set. For $M>0$ and $\d>0$, we say that a $d-$dimensional set $E$ is $(M,\d)$-quasiminimal on $U$ ($E\in QM(U,M,\d)$ for short) if $E$ is relatively closed in $U$ and for all $\d-$deformation $\phi_t$ on $U$ we have
\be H^d(E\cap W_1)\le M H^d(\phi_1(E\cap W_1),\ee
where a $\d$-deformation on $U$ is a family of maps $\{\phi_t\}_{t\in[0,1]}$ such that

1) $\phi_0=Id$ and $\phi_1$ is Lipschitz;

2) The map $(t,x)\mapsto \phi_t(x)$ is continuous from $[0,1]\times U$ to $U$;

3) If we denote
\be W_t=\{x\in U:\phi_t(x)\ne x\}\mbox{ and }W_\phi=\cup_{t\in[0.1]}W_t\cup\phi_t(W_t),\ee
then $W_\phi$ is relatively compact in $U$;

4) Diam$(W_\phi)<\d$.
 
\end{defn}

We want to prove that $F_k''$ is $(M',\d)$-quasiminimal in $\Omega_\e$, for some $M'$ and $\d$ that do not depend on $k$. 

First of all, by the construction of $F_k''$, and Lemma 5.2.2 of \cite{Fv}, we already know that $F_k''\in QM({ G}^\circ, A_0,+\infty)$ for some $A_0$ that does not depend on $k$. 

So let us take $\d=\min\{r,\frac{\e}{10000}\}$. Recall that $r>0$ is such that $ G$ is a Lipschitz neighborhood retract of $ G_r$. Denote by $\varphi$ such a neighborhood retract, with Lipschitz constant $C$.

Take any $\d-$deformation $\phi_t$ of $\Omega_\e$.  Set $V_1=W_1\cup\phi_1(W_1)$.

If $V_1$ is contained in $G^\circ$, then since $F_k''\in QM({ G}^\circ, A_0,+\infty)$, we have the desired estimate;

If $V_1$ does not meet $F_k''$, then also we have nothing to prove;

So the rest is to look at the case that $V_1$ is not contained in $ G^\circ$, but it meets $F_k''$. In this case $V_1\subset G_r$, and its intersection with $ G$ is contained in the region that is the union of dyadic cubes in $\cal D$.

We want to compare the measures of $W_1\cap F_k''$ and $\phi(W_1\cap F_k'')$. First we send the part of $V_1\bs G$ back to $G$ where we control things well, we have
\be H^d(\varphi\circ\phi(W_1\cap F_k''))\le C^dH^d(\phi(W_1\cap F_k'')).\ee
Notice that $W_1\subset G$, hence $\varphi$ does not move it.

Denote by $\cal D'$ all the dyadic cubes in $\cal D$ which touch a cube in $\cal D$ that touches $\varphi(V_1)$. Set $G_1=\cup_{Q\in\cal D'}Q$. By Lemma 4.3.2 of \cite{Fv}, we can find a Federer-Fleming projection $\pi$ from $\varphi\circ\phi(F_k''\cap G_1)$ to the $d-$skeleton of the union of cubes in $\cal D'$, which is $L'$-Lipschitz. Since $\varphi\circ\phi(W_1\cap F_k'')\subset \varphi\circ\phi(F_k''\cap G_1)$, we have
\be H^d(\pi\circ\varphi\circ\phi(W_1\cap F_k''))\le L'^dH^d(\varphi\circ\phi(W_1\cap F_k'')).\ee

But by construction of $F_k''$, the part $V=(F_k''\cap G_1)\bs W_1=\varphi\circ\phi [(F_k''\cap G_1)\bs W_1]\subset \varphi\circ\phi(F_k''\cap G_1)$ is part of $d-$skeleton of the unions of cubes in $\cal D'$, hence $\pi$ does not move it (c.f. \cite{DS00}, Proposition 3.1). Thus $V\cup\pi\circ\varphi\circ\phi(W_1\cap F_k'')$ is a deformation of $F_k''$ in $G_1$, and lives on the union of $d-$skeleton of $\cal D$. So by the $J_{g_k}$ minimality of $F_k''$, we have
\be \begin{split}J_{g_k}(V)+J_{g_k}(\pi\circ\varphi\circ\phi(W_1\cap F_k''))&\ge J_{g_k}(V\cup\pi\circ\varphi\circ\phi(W_1\cap F_k''))\\
&\ge J_{g_k}(F_k''\cap G_1)=J_{g_k}(W_1\cap F_k'')+J_{g_k}(V),\end{split}\ee
which gives 
\be J_{g_k}(\pi\circ\varphi\circ\phi(W_1\cap F_k''))\ge J_{g_k}(W_1\cap F_k'').\ee

Notice that $1\le g_k\le 10L^dM$, hence we have
\be (10L^dM)H^d(\pi\circ\varphi\circ\phi(W_1\cap F_k''))\ge H^d(W_1\cap F_k'').\ee
Combining with (4.14) and (4.13) we get
\be H^d(W_1\cap F_k'')\le (10L^dM)L'^dC^dH^d(\phi(W_1\cap F_k'')).\ee
Hence we have proved that $F_k''$ is $( M',\d)$ quasiminimal in $\Omega_\e$, where $M'=\max\{A_0,(10L^dM)L'^dC^d\}$.

Now we extract a converging subsequence of $F_k''$, still denoted by $F_k''$, and denote its limit by $F_0\in\overline U$.

By the uniform quasiminimality of $F_k''$, we have the lower semi continuity of $H^d$ (c.f.\cite{GD03} Theorem 3.4), and hence for all continuous function $f$, $J_f(F_0)\le\liminf_{k\to\infty}J_f(F_k')$. In particular for each $n$, 
\be J_{g_n}(F_0)\le\liminf_{k\to\infty}J_{g_n}(F_k'').\ee

Now for each $n$ fixed, when $k$ is large enough, we have $g_n\le g_k$, and hence $J_{g_n}(F_k'')\le J_{g_k}(F_k'')$, thus we have, by (4.19),
\be J_{g_n}(F_0)\le\liminf_{k\to\infty}J_{g_k}(F_k''),\ee
and by (4.9),
\be J_{g_n}(F_0)\le\liminf_{k\to\infty}(1+2^{-k})J_{g_k}(F_k)=\liminf_{k\to\infty}J_{g_k}(F_k)\le\liminf_{k\to\infty}J_h(F_k),\ee
since $F_k\subset\overline\Omega$.
Now since $g_n\uparrow h'$, the monotonous converging theorem gives
\be J_{h'}(F_0)=\lim_{n\to\infty}J_{g_n}(F_0)\le \liminf_{k\to\infty}J_h(F_k)=m.\ee

We still have to prove that 
\be F_0\in\mathfrak F.\ee 

First we prove that $F_0$ satisfies the topological condition with respect to $w_j,j\in J$. 

Since each $F_k''$ is a deformation of $F_k$ which is in $\mathfrak F$, hence by Proposition 3.6, $F_k''$ satisfy also that for $j\in J$, $w_j$ is not zero in $H_{n-d-1}(\R^n\bs F_k'')$. Now $F_0$ is the limit of $F_k''$. If there exists a $j\in J$ such that $w_j$ is zero in $H_{n-d-1}(\R^n\bs F_0)$, then there exists a smooth simplicial $n-d$-chain $\Gamma\in\R^n\bs F_0$ such that $\partial\Gamma=w_j$. Since $\R^n\bs F_0$ is open, and the support $|\Gamma|$ of $\Gamma$ is compact, there exists a neighborhood $V$ of $|\Gamma|$ such that $\overline V\cap F_0=\emptyset$. Then since $F_k''\to F_0$, and $\R^n\bs\overline V$ is open, there exists $N>0$ such that for all $k>N$, $F_k''\cap\overline V=\emptyset$, i.e. $\Gamma\subset\R^n\bs F_k$, too. In this case, $w_j$ is also zero in $H_{n-d-1}(\R^n\bs F_k''),$ for $k>N$. This is impossible.

So $F_0$ is such that $w_j$ is non zero in $H_{n-d-1}(\R^n\bs F_0)$ for all $j\in J$.

The last thing to prove is that 
$F_0\subset\Omega$, or equivalently, $H^d(F_0\bs\Omega)=0.$

Suppose this is not true, that is, $H^d(F_0\bs\Omega)=\a>0$. Notice that $F_0$ is the limit of $F_k''$, which are all contained in $G$, so $F_0$ is contained in $G$, and hence contained in $\Omega_\frac\e2$. Thus we can apply the $L-$neighborhood retract $\pi_\e$ (of $\Omega_\e$ to $\Omega$) to $F_0$, and $\pi_\e(F_0)$ is a deformation of $F_0$ in $\Omega_\e$, hence keeps all the $w_i$ non zero. Moreover, $\pi_\e(F_0)\subset \Omega$. Thus $\pi_\e(F_0)\in\mathfrak F.$

Let us calculate $J_h(\pi_\e(F_0))=J_{h'}(\pi_\e(F_0))$. We have
\be \begin{split}J_{h'}(\pi_\e(F_0))-J_{h'}(F_0)&=J_{h'}(\pi_\e(F_0\bs\Omega))-J_{h'}(F_0\bs\Omega)\\
&=\int_{\pi_\e(F_0\bs\Omega)}h'dH^d-\int_{F_0\bs\Omega}h'dH^d\\
&=\int_{\pi_\e(F_0\bs\Omega)}hdH^d-\int_{F_0\bs\Omega}10L^dMdH^d\\
&\le MH^d(\pi_\e(F_0\bs\Omega))-10L^dM H^d(F_0\bs\Omega)\\
&\le ML^dH^d(F_0\bs\Omega)-10L^dM H^d(F_0\bs\Omega)\\
&\le-9L^dMH^d(F_0\bs\Omega)=-9L^dM\a.
\end{split}\ee
That is
\be J_h(\pi_\e(F_0))=\le J_{h'}(F_0)-9L^dM\a\le m-9L^dM\a<m.\ee
This contradicts the fact that $m=\inf_{F\in\mathfrak F}J_h(F).$

Hence $F_0\subset\Omega$.

Thus we have (4.23), and hence $J_h(F_0)=J_{h'}(F_0)\le m$. But $F_0\in\mathfrak F$ gives the equality, thus the proof of Theorem 4.2 is finished.\qed

Next we are going to prove an existence theorem for a minimal topological competitor in a ball. 

\begin{thm}Let $E\subset \R^n$ be a closed set. $B\subset \R^n$ is an open ball. Let $\{w_j\}_{j\in J}$ be a family of smooth $n-d-1$-surfaces in $\R^n\bs(\overline B\cup E)$, which are non-zero in $H_{n-d-1}(\R^n\bs E)$. Set
\be \mathfrak F=\{F\subset \R^n,F\bs B=E\bs B\mbox{ and for all }j\in J, w_j\mbox{ is non-zero in }H_{n-d-1}(\R^n\bs F)\}.\ee

Then there exists $F_0\in\mathfrak F$ such that
\be H^d(F_0\cap B)=\inf\{H^d(F\cap B);F\in\mathfrak F\}.\ee
\end{thm}

\nd The idea is similar to the proof of Theorem 4.2. But we are going to use the uniform convexity of the ball $B$.

So set $m=\inf\{H^d(F);F\in\mathfrak F\}$. Fix a minimizing sequence $\{F_k\}_{k\in\N}\subset\mathfrak F$, that is
\be \lim_{k\to\infty}H^d(F_k\cap B)=m.\ee

Set $K=E\cap\partial B$, and set $U=\R^n\bs K$. Denote by $\mathfrak G$ the class of all deformations of one of the $F_k\cap B, k\ge 1$ in $U$. Then by Theorem 6.1.7 of \cite{Fv}, there exists a set $F_\infty\subset U$, which is the limit of some sequence $F_k'\in \mathfrak G$, such that 
\be H^d(F_\infty)\le\lim_{k\to\infty}H^d(F_k')=\inf_{F\in\mathfrak G}H^d(F)\le \inf_k H^d(F_k\cap B)=m.\ee
Moreover, $F_\infty$ is minimal in $U$.

We want to show that $F_\infty$ is contained in the convex hull $C_K$ of $K$. We need the following lemma.

\begin{lem}
Let $C\subset \R^n$ be a closed convex set with non-empty interior. Then for all $\e>0$, there exists $\d>0$ and a 1-Lipschitz retraction $f_\e$ of $\R^n$ onto $C$ such that $f_\e$ is $1-\d$-Lipschitz on $\R^n\bs B(C,\e)$.
\end{lem}

\nd

We define $f_C:\R^n\to\R$, for any $x\in\R^n\bs \{0\}$, $f_C(x)=\inf\{t:t^{-1}x\in C\}$, and $f_C(0)=0$. We have immediately that for any $\lambda>0$ and $x\in \R^n$, $f_C(\lambda x)=\lambda f_C(x)$. (If $C$ is symmetric, then $f_C$ is the norm on $\R^n$ whose unit ball is $C$. )

Then $f_C$ is convex. In fact, for all $x,y\in\R^n$, we have $\frac{x}{f_C(x)},\frac{y}{f_C(y)}\in C$ by definition. Thus for all $\a,\beta>0,\a+\beta=1$, we have
\be\begin{split} &\frac{\a x+\beta y}{\a f_C(x)+\beta f_C(y)}=\frac{x}{f_C(x)}\frac{f_C(x)\a}{\a f_C(x)+\beta f_C(y)}+\frac{y}{f_C(y)}\frac{\beta f_C(y)}{\a f_C(x)+\beta f_C(y)}\\
&=\frac{x}{f_C(x)}\frac{f_C(x)\a}{\a f_C(x)+\beta f_C(y)}+\frac{y}{f_C(y)}(1-\frac{f_C(x)\a}{\a f_C(x)+\beta f_C(y)})\in C,\end{split}\ee
therefore $f_C(\a x+\beta y)\le \a f_C(x)+\beta f_C(y).$

Since $C^\circ\ne\emptyset$, there exists $A>1$ such that
\be A^{-1}f_C\le||\cdot||\le Af_C\ee
where $||\cdot||$ denotes the Euclidean norm.

Now for all $a\ge 0$, set $f_a=f_C+a||\cdot||$ and $C_{a,b}$ the closed set $\{x\in\R^n:f_a(x)\le b\}$. Then $f_a$ is convex since $f_C$ and $||\cdot||$ are. Thus $C_{a,b}$ is also convex.  Notice that for all $x\in\R^n\bs\{0\}$, $f_a(x)$ is a strictly increasing continuous function of $a$, $f_0=f_C$, and 
\be C_{a,b}\supset C_{a',b},\ C_{a,b}\subset C_{a,b'}\mbox{ for all }a<a',b<b',\ee
\be \bigcap_{a_n\to a-}C_{a_n,b}=\bigcap_{b_n\to b+}C_{a,b_n}=C_{a,b};\bigcup_{a_n\to a+}C_{a_n,b}=\bigcup_{b_n\to b-}C_{a,b_n}=C_{a,b}^\circ.\ee

Now since the $f_a$ contains a part of Euclidean norm, which is uniformly convex, it is easy to verify that
\be\begin{split}\mbox{for all }&a,b>0, \mbox{ there exists a constant }M(a,b,A)\mbox{ such that}\\
&\mbox{ for each }x,y\in \partial C_{a,b}\mbox{ with }\a_{x,y}<\frac\pi2,
 B(\frac{x+y}{2}, M(a,b,A)||x-y||^2)\subset C_{a,b},\end{split}\ee
 with $M(a,b,A)>0$, where $\a_{x,y}<\pi$ denotes the angle between $\vec{Ox}$ and $\vec{Oy}$ for $x,y\ne 0$, and $B(\frac{x+y}{2}, M(a,b,A)||x-y||^2)$ denotes the euclidean ball centered at $\frac{x+y}{2}$ with radius $M(a,b,A)||x-y||^2$.

Now for all $\e>0$, let $w,v$ be two points in $\R^n\bs B(C_{a,b},\e)$ such that $\pi_{a,b}(w)=x,\pi_{a,b}(v)=y$, where $\pi_{a,b}$ denotes the shortest distance projection on the convex set $C_{a,b}$. We claim that the angle $\beta_1\in[0,\frac\pi 2]$ between $\vec{xw}$ and $\vec{yx}$ is smaller than $\arctan\frac{1}{2M(a,b,A)||x-y||}.$ (See the picture 4-1 below). In fact, denote by $P$ the plane containing $x,y$ and $w$, denote by $z\in P$ the point such that $[z,\frac{x+y}{2}]\perp [x,y]$ and $||z-\frac{x+y}{2}||=M(a,b,A)||x-y||^2$. Then $z\in C_{a,b}$, $[x,z]\in C_{a,b}$, and
\be\tan\angle zxy=2M(a,b,A)||x-y||.\ee

Then if $\beta_1>\arctan\frac{1}{2M(a,b,A)||x-y||}$, we have $\angle wxz<\frac\pi2$. Denote by $s$ the projection of $w$ on $L'$, the line passing through $x$ and $z$. Then $s$ is between $x$ and $z$, or $z$ is between $x$ and $s$. In both cases $(x,z)\cap (x,s)\ne \emptyset$. Take $x'\in  (x,z)\cap (x,s)\subset C_{a,b}$, then $x'\in C_{a,b}$, and $\angle wxz<\angle wx'z$. As a result
\be ||w-x||=||w-s||/\sin\angle wxz>||w-s||/\sin\angle wx'z=||w-x'||,\ee
which contradicts the fact that $x$ is the shortest distance projection of $w$ on $C_{a,b}.$

\noindent\centerline{\includegraphics[width=0.8\textwidth]{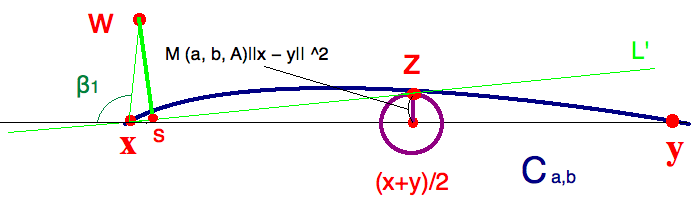}}
\nopagebreak[4]
\centerline{4-1}

Similarly we can prove that if $\beta_2$ denotes the angle between $\vec{yv}$ and $\vec{xy}$, then $\beta_2\le \arctan\frac{1}{2M(a,b,A)||x-y||}.$

Denote by $L$ the line passing through $x$ and $y$, $Q$ the plane perpendicular to $L$, $p_l$ and $p_Q$ orthogonal projections on them. Then
\be ||w-v||\ge ||p_L(w)-p_L(v)||=||w-x||\cos\beta_1+||x-y||+||v-y||\cos\beta_2.\ee
But $w,v\in\R^n\bs B(C_{a,b},\e)$, hence we have $||w-x||>\e,||v-y||>\e$, therefore
\be \begin{split}||w-v||&\ge ||x-y||+2\e\cos\arctan\frac{1}{2M(a,b,A)||x-y||}\\
&=(1+\e C(a,b,A))||x-y||,\end{split}\ee
with $C(a,b,A)>0$.

Notice that (4.40) is true for all pairs of $x,y$ such that $\a_{x,y}<\frac\pi2$. Hence $\pi_{a,b}$ is locally $\frac{1}{1+\e C(a,b,A))}$-Lipschitz on $\R^n\bs B(C_{a,b},\e)$.

Return to the proof of the lemma. Fix an arbitrary $\e>0$. Then by (4.34) and (4.35), there exists $a,b>0$ such that
\be C\subset C_{a,b}\subset B(C_{a,b},\frac\e2)\subset B(C,\e).\ee

Now denote by $\pi_C$ the shortest distance projection on $C$. Denote by $f_\e=\pi_C\circ \pi_{a,b}$ for a pair of $a,b$ which satisfies (4.41). Then for proving the lemma, it is sufficient to prove that $\pi_{a,b}$ is locally $1-\d$-Lipschitz on $\R^n\bs B(C_{a,b},\frac\e2)$. Then by (4.40), we take $\d$ such that $1-\d=\frac{1}{1+\frac12 \e C(a,b,A))}$, and we obtain the conclusion.\qed

\begin{cor}
Let $E\subset \R^n\bs B(C,\e)$ be a rectifiable set and $f$ be as in the lemma, then
\be H^d(f(E))\le (1-\d)^dH^d(E).\ee
\end{cor}

Now we can prove that $F_\infty\subset C_K$. Recall that $C_K$ is the convex hull of $K$.

We will prove that 
\be\mbox{for all }\e>0,F_\infty\subset B(C_K,2\e).\ee

Suppose this is not true for a $\e>0$, i.e. $H^d(F_\infty\bs B(C_K,2\e))>0$. We apply Corollary 4.42 to $F_\infty\bs B(C_K,2\e)$ and the convex set $B(C_K,\e)$, and obtain that there exists a Lipschitz map $f_\e$ of $\R^n$ in $B(C_K,\e)$, such that
\be H^d(f_\e(F_\infty))<H^d(F_\infty),\ee
where $f_\e$ is as in Lemma 4.31. But $f_\e(F_\infty)$ is a deformation of $F_\infty$ in $U$. In fact $F_\infty$ is compact, hence $F_\infty$ is contained in a ball $B(0,R)$. Define a map $f:U\to U$, $f=f_\e$ in $B(0,R)$, and $f=id$ in $\R^n\bs B(0,2R)$. Then $f(F_\infty)=f_\e(F_\infty)$. But the set $\{x\in U: f(x)\ne x\}$ is compact in $U$, because $f_\e$ does not move points near $C_K$, which contains $K=\partial U$. Hence $f(F_\infty)$ is a deformation of $F_\infty$ in $U$, and hence $f_\e(F_\infty)$ is, too.

Thus if $H^d(F_\infty\bs B(C_K,2\e))>0$, then we can find one of its deformation in $U$ which decrease its measure, this contradicts the fact that $F_\infty$ is minimal in $U$.

Hence we have (4.44), and thus we have $F_\infty\subset\cap_{\e>0} B(C_K,2\e)=C_K$. This gives that $F_\infty\subset B$, and $\overline F_\infty\cap\partial B\subset C_K\cap\partial B=K$. But $\overline F_\infty\supset K$, hence $\overline F_\infty\cap\partial B=K$.

Set $F_0=F_\infty\cup (E\bs B)$. We want to prove that $F_0\in\mathfrak F$. We know that $F_\infty\subset B$, hence $F_0\bs B=E\bs B$. So we only have to verify the topological condition on $w_j,j\in J$.

So take any $j\in J$. Denote by $r=\mbox{dist}(w_j,\overline B)<\infty$, and $B'=B(B,\frac12 r)$ the concentric ball of $B$ with a radius $\frac12 r$ larger than that of $B$. Suppose that $w_j$ is zero in $H_{n-d-1}(\R^n\bs F_0)$. Then there exists a smooth simplicial $n-d-$chain $\Gamma\subset\R^n\bs F_0$ such that $\partial\Gamma=w_j$. Now since $F_0\cap\overline B=\lim_{k\to\infty} F_k'$, when $k$ is large, we have $F_k'\subset B'\bs\Gamma$. Since $F_k'\subset\mathfrak G$, there exists a $l\ge 1$ such that $F_k'$ is the deformation of $F_l\cap  B$ in $U$. But $F_k'\subset B'$, hence $F_k'$ is a deformation of $F_l\cap  B$ in $B'\bs K$. And hence $F_k'\cup (F_l\bs B)$ is a deformation of $F_l$ in $B'$. By Proposition 3.7, $F_k'\cup (F_l\bs B)$ is a topological competitor of $F_l$ in $B'$. But note that $F_k'\cup (F_l\bs B)$ does not touch $\Gamma$, hence $\partial \Gamma=w_j$ is zero in $H_{n-d-1}(\R^n\bs(F_k'\cup (F_l\bs B))$, and therefore $w_j$ is zero in $H_{n-d-1}(\R^n\bs F_l)$, which contradicts the fact that $F_l\in \mathfrak F$.

Hence $w_j$ is non-zero in $H_{n-d-1}(\R^n\bs F_0)$, for any $j\in J$. Thus we get $F_0\in\mathfrak F$. Now by (4.30),
\be H^d(F_0\cap B)=H^d(F_\infty)\le m=\inf_{F\in\mathfrak F}H^d(F\cap B),\ee
hence we have the equality, that is
\be  H^d(F_0\cap B)=\inf_{F\in\mathfrak F}H^d(F\cap B).\ee
Thus the proof of Theorem 4.26 is completed.\qed

%
%
%

\begin{rem} Theorem 4.26 is still true with the same proof if we replace the ball $B$ by any convex set  whose boundary points are all extremal points.
\end{rem}

As an application, we give the following corollary, for preparing for the next section.
\begin{cor}
Let $P,Q$ be two planes contained in $\R^4$, with $P\cap Q=\{0\}$. Denote by $B=B(0,1)$ the unit ball. Denote by $\mathfrak F$ all the topological competitors of $E=P\cup Q$ in $B$. Then there exists $F_0\in \mathfrak F$ such that 
\be H^2(F_0\cap B)=\inf_{F\in\mathfrak F}H^2(F\cap B).\ee
Moreover $F_0\cap\overline B$ is contained in the convex hull of $E\cap \overline B$. 
\end{cor}

\nd This is a direct corollary of Theorem 4.26, where we take $\{w_j\}_{j\in J}$ to be the family of all circles outside $B$ who are non-zero in $H_1(\R^n\bs (P\cup Q))$.\qed

\begin{rem}We could also get Corollary 4.49 by using Theorem 4.2. Denote by $C$ the convex hull of $(P\cup Q)\cap \overline B$. Intuitively we want to take $P\cup Q\cup C$ to be the $\Omega$ in Theorem 4.2, $h=1$, and $\{w_j\}_{j\in J}$ being all smooth closed curves in $\R^4\bs \Omega$ which are not zero in $H_1(\R^4\bs\Omega)$. 

The problem is $\Omega$ is not compact, and the class of $\{w_j\}_{j\in J}$ may be arbitrarily close to $\Omega.$

But it is not hard to prove that every $w_j,j\in J$, has a representative which is on the sphere $\partial B(0,2)$, and is $\frac{1}{100}$ far from $\Omega$. 

Then to deal with the compactness, first we add a point $\infty$ such that $\overline{\R^4}=\R^4\cup\{\infty\}$ is the one point compactification of $\R^4$, and $E=P\cup Q\cup\{\infty\}$ is compact in $\overline{\R^4}$. Then $\overline{\R^4}\bs E=\R^4\bs(P\cup Q)$. Take any point $p\in\R^4\bs B(0,10)$, then $H_1(\overline{\R^4}\bs E)\cong H_1((\overline{\R^4}\bs\{p\})\bs E)$. Next we take a homeomorphism $\varphi$ from $\overline{\R^4}$ to itself, with $\varphi|_{\overline B(0,5)}=id$, and $\varphi(p)=\infty$. Then since $w_j,j\in J$ are all in $\partial B(0,2)$, they are non-zero in $\R^4\bs (P\cup Q)$ if and only if they are non-zero in $\R^4\bs\varphi(E)$. Notice that for guaranteeing that all $w_j$ to be non-zero, any $F\in\mathfrak F$ should contain the part $P\cup Q\bs C$, because a hole will make a local non-zero represent zero. Hence what we want to minimize is the part of $F$ inside $C$ over all $F\in\mathfrak F$, which does not change by $\varphi$. This time $\varphi(E)$ is compact, we denote it by $\Omega$, and we can apply Theorem 4.2 to get the desired conclusion.
\end{rem}

\section{The union of two almost orthogonal planes is topologically minimal in $\R^4$}

This section is devoted to proving the topological minimality of the union of two almost orthogonal planes in $\R^4$. More precisely, we are going to prove the following theorem.

\begin{thm}[Topological minimality of the union of two almost orthogonal planes] There exists $0<\theta<\frac\pi2$, such that if $P^1$ and $P^2$ are two planes in $\R^4$ with characteristic angles $\a_1\le\a_2$ and $\a_1\ge\theta$, then their union $P^1\cup P^2$ is a topological minimal cone.
\end{thm}

Here the characteristic angles between two planes $P,Q$ are defined as follows : 

$\a_1=\min\{\mbox{angle between }v\mbox{ and }w;v\in P,w\in Q\mbox{ are unit vectors}\}$. We fix $v_1\in P$ and $w_1\in Q$ such that the angle between $v_1$ and $w_1$ is $\a_1$. Then $\a_2$ is defined as $\min\{\mbox{angle between }v\mbox{ and }w;v\in P,w\in Q\mbox{ are unit vectors}, v\perp v_1,w\perp w_1\}$.

Notice that $\a_1=\a_2=\frac\pi2$ means that the two planes are orthogonal.

The general idea is the same as that in \cite{2p}, which proved the Almgren minimality of almost orthogonal planes. So we will keep the main structure of the proof here, and just give proofs for the places where things are different, especially for the parts that concerning projections.

And also, the conclusion can be generalized to the union of several almost orthogonal $m-$planes, for any $m\ge 2$. Things are different (and sometimes more complicated) in some places, such as the uniqueness theorem, and the harmonic extension. But we are not going to discuss them here. Please refer to \cite{XY10} for details. 

\subsection{Some basic preliminaries and estimates for unit simple 2-vectors}

Denote by $\wedge_2(\R^4)$ the space of all 2-vectors in $\R^4$. Let $x,y$ be two vectors in $\R^4$, we denote by $x\wedge y\in\wedge_2(\R^4)$ their exterior product. And if $\{e_i\}_{1\le i\le 4}$ is a orthonormal basis, then $\{e_i\wedge e_j\}_{1\le i<j\le 4}$ forms a basis of $\wedge_2 (\R^4)$. We say that an element $v\in\wedge_2 (\R^4)$ is simple if it can be expressed as the exterior product of two vectors. 

The norm on $\wedge_2 (\R^4)$, denoted by $|\cdot |$, is defined by 
\be|\sum_{i<j}\lambda_{ij}e_i\wedge e_j|=\sum_{i<j}|\lambda_{ij}|^2.\ee
Under this norm $\wedge_2(\R^4)$, is a Hilbert space, and $\{e_i\wedge e_j\}_{1\le i<j\le 4}$ is an orthonormal basis. For all simple 2-vector $x\wedge y$, its norm is 
\be |x\wedge y|=||x||||y||\sin<x,y>,\ee
where $<x,y>\in[0,\pi]$ is the angle between the vectors $x$ and $y$, and $||\cdot ||$ denotes the Euclidean norm on $\R^4$. A unit simple 2-vector is a simple 2-vector of norm 1. Notice that $|\cdot|$ is generated by the scalar product $<,>$ defined on $\wedge_2(\R^4)$ as the following: for $\xi=\sum_{1\le i<j\le 4} a_{ij}e_i\wedge e_j, \zeta=\sum_{1\le i<j\le 4} b_{ij}e_i\wedge e_j$, 
\be <\xi,\zeta>=\sum_{1\le i<j\le 4} a_{ij}b_{ij}.\ee

One can easily verify that if two pairs of vectors $x,y$ and $x',y'$ generate the same 2-dimensional subspace of $\R^4$, then there exists $r\in\R\bs\{0\}$ such that $x\wedge y=rx'\wedge y'$.

Now given a unit simple 2-vector $\xi$, we can associate it to a 2-dimensional subspace $P(\xi)\in G(4,2)$, where $G(4,2)$ denote the set of all 2-dimensional subspaces of $\R^4$:
\be P(\xi)=\{v\in\R^4,v\wedge\xi=0\}.\ee 
In other words, $P(x\wedge y)$ is the subspace generated by $x$ and $y$.

(From time to time, when there is no ambiguity, we write also $P=x\wedge y$, where $P\in G(4,2)$ and the two unit vectors $x,y\in\R^4$ are such that $P=P(x\wedge y)$. In this case $x\wedge y$ represents a plane.)

For the side of linear maps, if $f$ is a linear map from $\R^4$ to $\R^4$, then we denote by $\wedge_2f$ (and sometimes $f$ if there is no ambiguity) the linear map from $\wedge_2(\R^4)$ to $\wedge_2(\R^4)$ such that
\be \wedge_2f(x\wedge y)=f(x)\wedge f(y).\ee

And for the side of $G(4,2)$ (the set of all planes, without considering orientations), for a unit simple 2-vector $\xi\in\wedge_2\R^4$, we have always $P(\xi)=P(-\xi)$, so that we can define $|f(\cdot)|:G(4,2)\to\R^+\cup\{0\}$ by
\be |f(P(\xi))|=|\wedge_2f(\xi)|.\ee
One can easily verify that the value of  $|f(P(\xi))|$ does not depend on the choice of $\xi$ that generates $P$. So $|f(\cdot)|$ is well defined.

\subsection{First results concerning unions of planes}

\begin{lem}[c.f.\cite{2p} lemmas 2.15 and 2.18] For every unit simple 2-vector $\xi\in\wedge_2\R^4$,
\be |p^1(\xi)|+|p^2(\xi)|\le 1,\ee
and if we denote by $\Xi$ the subset of all unit simple 2-vectors that satisfies the equality, then for any two unit vectors $x,y$, $x\wedge y\in\Xi$ if and only if there exist
 $\a\in[0,\frac\pi2]$, $v_i,u_i,i=1,2$ four unit vectors, $v_i\in P^1,\ u_i\in P^2$, $v_1\perp v_2,u_1\perp u_2$, such that
\be x=\cos\a v_1+\sin\a u_1\mbox{ and }y=\cos\a v_2+\sin\a u_2.\ee
\end{lem}

\begin{lem}[c.f. \cite{2p} Proposition 2.20] Let $0\le\alpha_1\le \alpha_2\le\frac\pi2$, and let $P^1,P^2\subset\R^4$ be two planes with characteristic angles $\alpha_1\le \alpha_2$. Denote by $p^i$ the orthogonal projection on $P^i$, $i=1,2$. Then for all unit simple 2-vector $\zeta\in \bigwedge_2\R^4$, its projections on these two planes satisfy:
\be |p^1\zeta|+|p^2\zeta|\le 1+2\cos\alpha_1.\ee
\end{lem}

\begin{pro}
 Let $P^1$ and $P^2$ be two orthogonal planes in $\R^4$. Then $P_0:=P^1\cup P^2$ is topologically minimal.
\end{pro}

\nd

Since $P_0$ is a cone, so we just have to prove that $P_0$ is locally topologically minimal in the unit ball $B$.

The idea is just to prove that the projections to $P^i,i=1,2$ of $E$ should be surjective in $\overline B$, and then we are satisfied because of the Wirtinger's inequality.

Denote by $p^i,i=1,2$ the orthogonal projection on $P^i$. Suppose for example that $p^1$ is not surjective in $\overline B$, then there exists $x\in P^1\cap\overline B$ such that ${p^1}^{-1}(x)\cap E\cap\overline B=\emptyset$. Since $E$ is a competitor of $P_0$ in $B$, hence $E\cap\partial B=P_0\cap\partial B$, hence $x\in B\bs\{0\}$. But in this case, ${p^1}^{-1}(x)\cap (E\bs B)={p^1}^{-1}(x)\cap (P_0\bs B)=\emptyset$. So we have
\be {p^1}^{-1}(x)\cap E=\emptyset.\ee

Notice that ${p^1}^{-1}(x)=x+P^2$. Denote by $s=(x+P^2)\cap\partial B(0,2)$ a circle in $\R^4\bs (B\cup P_0)$. Then $s$ is contractible in $x+P^2$. But $E\cap (x+P^2)=\emptyset$, so $s$ is contractible in $\R^4\bs E$. In other words, $s$ is zero in $H_1(\R^4\bs E)$. But obviously $s$ is not zero in $H_1(\R^4\bs P_0)$ (it is an generator of the group $H_1(\R^4\bs P^1)$). This contradicts the fact that $E$ is a competitor of $P_0$ with respect to $B$.

So \be p^i(E)\supset P^i\cap B.\ee

Now we cite the following lemma, which will help us to finish the prove of Proposition 5.13.

\begin{lem}[c.f.\cite{2p}, Lemma 2.31]\label{projection}Let $P^1, P^2$ be two planes in $\R^4$, let $F\subset\R^4$ be a 2-rectifiable set. Denote by $p^i$ the projection on $P^i$. If $\lambda$ is such that for almost all
$x\in F$, the approximate tangent plane of $F$ at $x$ $T_xF\in G(4,2)$ satisfies that
\be |p^1(T_xF)|+|p^2(T_xF)|\le \lambda,\ee 
then we have 
\be H^2(p^1(F))+H^2(p^2(F))\le \lambda H^2(F).\ee
\end{lem} 

On combining Lemmas 5.8 and 5.16, we have that for our set $E$,
\be H^2(E)\ge H^2(p^1(E))+H^2(p^2(E))\ge H^2(P^1\cap B)+H^2(P^2\cap B)=H^2(P_0\cap B).\ee

This is true for any topological competitor $E$ of $P_0$ in $B$. Hence $P_0$ is topological minimal.\qed

\begin{cor} [uniqueness of $P_0$] Let $P_0=P_0^1\cup_\perp P_0^2$, and denote by $p_0^i$ the orthogonal projection on $P_0^i,i=1,2$. Let $E\subset \overline B(0,1)$ be a 2-dimensional closed reduced set which is topologically minimal in $B(0,1)\subset\ \R^4$ and satisfies:
\be p_0^i(E\cap \overline B(0,1))\supset P_0^i\cap \overline B(0,1);\ee
\be E\cap \partial B(0,1)=P_0\cap \partial B(0,1);\ee
\be \begin{split}&H^2(E\cap B(0,1))\le 2\pi\\
\mbox{ (or equivalently }&H^2(E\cap B(0,1))=2\pi\mbox{ by  Lemmas 5.8, 5.16 and (5.21))}.\end{split}\ee

Then $E=P_0\cap \overline B(0,1)$.

\end{cor}

\nd This is a direct corollary of Proposition 5.13, the uniqueness theorem 3.1 in \cite{2p}, and the fact that all topological minimal sets are Almgren minimal.\qed.

\subsection{A converging sequence of topological minimal competitors}

We begin to prove the promised theorem 5.1 at the beginning of this section. We prove it by contradiction. So suppose that the theorem is not true. That is, for any $k\in\N$, there exists two planes $P_k^1$ and $P_k^2$, with characteristic angles $\frac\pi2>\a_k^2\ge\a_k^1>\frac\pi2-\frac1k$, such that their union $P_k:=P_k^1\cup P_k^2$ is not topological minimal. Recall also that $P_0=P_0^1\cup P_0^2$ is the orthogonal union of two planes. Chose an orthonormal basis $\{e_i\}_{1\le i\le 4}$ of $\R^4$ such that $P_0^1=P(e_1\wedge e_2)$ and $P_0^2=e_3\wedge e_4$. After necessary rotations, we suppose also that all the $P_k^1,k\ge 0$ are the same, and $P_k^2=P((\cos\a_k^1 e_1+\sin\a_k^1 e_3)\wedge(\cos\a_k^2 e_2+\sin\a_k^2 e_4)$.

Then by the definition of topological minimal sets, and the fact that $P_k$ is a cone, we can find a topological competitor $E$ of $P_k$ with respect to the unit ball $B(0,1)$, such that
\be H^2(E\cap B(0,1))<H^2(P_k\cap B(0,1))=2\pi.\ee

Now by Corollary 4.49 in the previous section, there exists a competitor $F_k$ of $P_k$ with respect to $B=B(0,1)$, which is locally topologically minimal in $B$, and $E_k\cap\overline B$ is contained in the convex hull $C_k$ of $P_k\cap\overline B$. Denote by $E_k=F_k\cap\overline B$ the part of $F_k$ inside $\overline B$. The minimality of $E_k$ gives that 
\be H^2(E_k)<H^2(P_k\cap B(0,1))=2\pi.\ee

Now since $\overline B$ is compact, we can extract a converging subsequence of $\{E_k\}$, denoted still by $\{E_k\}$ for short. Denote by $E_\infty$ their limit. Then $E_\infty$ is contained in $\cap_n\cup_{k>n} C_k$, such that $E_\infty\cap \partial B\subset (\cap_n\cup_{k>n} C_k)\cap\partial B=P_0\cap\partial B$. On the other hand $E_\infty\cap \partial B\supset\lim_{k\to\infty} (P_k\cap\partial B)=P_0\cap\partial B$. Hence 
\be E_\infty\cap\partial B=P_0\cap\partial B.\ee

We want to use the uniqueness theorem 5.20, to prove that $E_\infty$ is in fact $P_0$. So we have to check all the conditions.

1) First let us check that $E_\infty$ is a topological competitor of $P_0$ in $B$. It is easy to see that $E_\infty\bs B=P_0\bs B$, in particular because $E_\infty\cap \partial B=P_0\cap\partial B$. Now if $E_\infty$ is not a competitor, then there exists a circle $s\in \R^4\bs (\overline B\cup P_0)$ such that $s$ is non zero in $H_1(\R^4\bs P_0)$, but is zero in $H_1(\R^4\bs E_\infty)$. That means, there exists a smooth simplicial 2-chain $\Gamma\subset\R^4\bs E_\infty$ such that $\partial \Gamma=s$. Then since $P_k\to P_0$ in any compact sets, we have for $k$ large, $s\in\R^4\bs P_k$ and $s$ is non zero in $H_1(\R^4\bs P_k)$. On the other hand since $E_k\to E_\infty$, when $k$ is large, $E_k\cap \Gamma=\emptyset$. Thus when $k$ is large enough, $s$ is zero in $\R^4\bs E_k$, which contradicts the fact that $E_k$ is a topological competitor of $P_k$.

2) Now since $E_\infty$ is a topological competitor of $P_0$ in $B$, we have automatically (5.21) by (5.15) The condition (5.22) has already been checked in (5.26). For (5.23), we know that $E_k$ is a sequence of Al-minimal sets, which are uniformly concentrated, so we have the lower-semicontinuity $H^2(E_\infty)\le\liminf_{k\to\infty}H^2(E_k)$. Thus (5.23) follows by (5.25).

So the rest is to prove that $E_\infty$ is topologically minimal. 

However, since $E_\infty$ is a competitor of $P_0$ in $B$, all competitors of $E_\infty$ in $B$ is automatically competitors of $P_0$ in $B$ (just by definition). But $H^2(E_\infty\cap\overline B)=2\pi$ implies already that it minimizes the Hausdorff measure among all competitors of $P_0$ in $B$, hence of course it minimizes the Hausdorff measure among all its competitors in $B$, hence it is topological minimal in $B$.

Thus we have checked all the conditions for the uniqueness theorem, and hence 
\be E_\infty=P_0.\ee

Thus we have a sequence of closed sets $E_k$, each $E_k$ is a minimal topological competitor of $P_k$ in $B$, and $E_k$ converges to $P_0$.

\subsection{A stopping time argument}

We will continue our argument, by cutting each $E_k$ into two pieces. One piece is inside a small ball near the origin, where something complicated happens there, and we can only estimate its measure by projection argument; the other piece is outside the small ball, where $E_k$ is very near $P_k$, and by the regularity of minimal sets near planes, $E_k$ is composed of two $C^1$ graphs on $P_k^1$ and $P_k^2$, where we will estimate their measures, which is equivalent to the Dirichlet energy, by harmonic extensions. So the first step is to find this small ball, with the critical radius, by a stopping time argument.

For each $k$ and $i=1,2$, denote by 
\be C_k^i(x,r)={p_k^i}^{-1}(B(0,r)\cap P_k^i)+x,\ee
and
\be D_k(x,r)=C_k^1(x,r)\cap C_k^2(x,r).\ee
Notice that $D_k(x,r)\supset B(0,1)$ and $D_k(0,1)\cap P_k=B(0,1)\cap P_k$.

We say that two sets $E,F$ are $\e r$ near each other in an open set $U$ if
\be d_{r,U}(E,F)<\e,\ee 
where
\be d_{r,U}(E,F)=\frac 1r\max\{\sup\{d(y,F):y\in E\cap U\},\sup\{d(y,E):y\in F\cap U\}\}.\ee
We define also
\be\begin{split}
    &d_{x,r}^k(E,F)=d_{r,D_k(x,r)}(E,F)\\
&=\frac 1r\max\{\sup\{d(y,F):y\in E\cap D_k(x,r)\},\sup\{d(y,E):y\in F\cap D_k(x,r)\}\}.
   \end{split}\ee

\begin{rem}
 Observe that $d_{r,U}(E,F)\ne\frac 1rd_H(E\cap U,F\cap U)$. For example, take $U=B(0,1)$, set $E_n=\partial B(0,1-\frac 1n)$, and $F_n=\partial B(0,1+\frac 1n)$, then we have \be d_{1,U}(E_n,F_n)\to 0, d_H(E_n\cap U,F_n\cap U)=d_H(E_n\cap U,\emptyset)=\infty.\ee
\end{rem}

Now we start our stopping time argument. We fix a $\e$ small and a $k$ large, and we set $s_i=2^{-i}$ for $i\ge 0$. Denote by $D(x,r)=D_k(x,r), d_{x,r}=d^k_{x,r}$ for short. Then we proceed as follows.

Step 1: Denote by $q_0=q_1=0$, then in $D(q_0,s_0)$, the set $E_k$ is $\e s_0$ near $P_k+q_1$ when $k$ is large, because $E_k\to P_0$ and $P_k\to P_0$ implies that $d_{0,1}(E_k,P_k)\to 0$.

Step 2: If in $D(q_1,s_1)$, there is no point $q\in\R^4$ such that $E_k$ is $\e s_1$ near $P_k+q$, we stop here; otherwise, there exists a point $q_2$ such that $E_k$ is $\e s_1$ near $P_k+q_2$ in $D(q_1,s_1)$. Here we ask $\e$ to be small enough (say, $\e<\frac {1}{100}$) such that such a $q_2$ is automatically in $D(q_,\frac12 s_1)$, by the conclusion of the step 1. Then in $D(q_1,s_1)$ we have simultaneously  
\be d_{q_1,s_1}(E_k,P_k+q_1)\le s_1^{-1}d_{q_0,s_0}(E_k,P_k+q_1)\le 2\e;\ d_{q_1,s_1}(E_k, P_k+q_2)\le\e.\ee
This implies that $d_{q_1,\frac12 s_1}(P_k+q_1,P_k+q_2)\le 12\e$ when $\e$ is small. And hence $d(q_1,q_2)\le 6\e$.

Now we are going to define our iteration process. Notice that this process depends on $\e$, hence we also call it the $\e-$process.

Suppose that $\{q_i\}$ are defined for all $i\le n$, with 
\be d_{q_1,q_{i+1}}\le 12 s_i\e=12\times 2^{-i}\e\ee
for $0\le i\le n-1$, and hence
\be d_{q_i,q_j}\le 24\e_{\min(i,j)}=2^{-\min(i,j)}\times 24\e\ee
for $0\le i,j\le n$, and that for all $i\le n-1$, $E_k$ is $\e s_i$ near $P_k+q_{i+1}$ in $D(q_i,s_i)$. We say in this case that the process does not stop at step $n$. Then

Step $n+1$: We look inside $D(q_n,s_n)$.

If $E_k$ is not $\e s_n$ near any $P_k+q$ in this ``ball'' of radius $s_n$, we stop. In this case, since $d(q_{n-1},q_n)\le 12\e s_{n-1}$, we have $D(q_n,2s_n(1-12\e))=D(q_n,s_{n-1}(1-12\e))\subset D(q_{n-1},s_{n-1})$, and hence
\be\begin{split}
    d_{q_n,2s_n(1-12\e)}(P_k+q_n,E_k)&\le (1-12\e)^{-1}d_{q_{n-1},s_{n-1}}(P_k+q,E_k)\\
&\le\frac{\e}{1-12\e}.
   \end{split}
\ee
Moreover
\be d(q_n,0)=d(q_n,q_1)\le 2^{-\min(1,n)}\times 24\e=12\e.\ee

Otherwise, we can find a $q_{n+1}\in\R^4$ such that $E_k$ is still $\e s_n$ near $P_k+q_{n+1}$ in $D(q_n,s_n)$, then since $\e$ is small, as before we have $d(q_{n+1},q_n)\le 12\e s_n$, and for $i\le n-1$,
\be d_{q_i,q_n}\le\sum_{j=i}^n 12\times 2^{-j}\e\le 2^{-\min(i,n)}\times 24\e.\ee
Thus we get our $q_{n+1}$, and say that the process does not stop at step $n+1$.

We will see in the next subsection, that the process has to stop at a finite step. And for each $k$, if the process stop at step $n$, we define $o_k=q_n,r_k=s_n$. Then $D_k(o_k,r_k)$ is the critical ball that we look for, because inside the small ball, by definition we know that $E_k$ is $\e s_n$ far from any translation of $P_k$, but outside it, things are near. We also have, by (5.39), $d(o_k,0)\le 12\e$, hence near the origin.

\subsection{Regularity and projection properties of $E_k$}

\begin{pro}\label{reg}
 There exists $\e_0\in(0,\frac{1}{100})$, such that for any $\e<\e_0$ fixed and for $k$ large, if our $\e-$process does not stop before the step $n$, then

(1) The set $E_k\cap (D_k(0,\frac{39}{40})\bs D_k(q_n,\frac{1}{10}s_n))$ is composed of two disjoint pieces $G^i,i=1,2$, such that 
\be G^i\mbox{ is the graph of a }C^1\mbox{ map }g^i:D_k(0,\frac{39}{40})\bs D_k(q_n,\frac{1}{10}s_n)\cap P_k^i\to {P_k^i}^\perp\ee
with
\be ||\nabla g^i||_\infty<1;\ee

(2) For each $t\in[\frac{1}{10}s_n,s_n],$
\be E_k\cap (D_k(0,1)\bs D_k(q_n,t))=G_t^1\cup G_t^2,\ee
where $G_t^1$ and $G_t^2$ do not meet. Moreover
\be P_k^i\cap (D_k(0,1)\bs C_k^i(q_n,t))\subset p_k^i(G_t^i)\mbox{ for }i=1,2,\ee
where $p_k^i$ is the orthogonal projection on $P_k^i,i=1,2$;

(3) The projections $p_k^i:E_k\cap \overline D_k(q_n,t)\to P_k^i\cap \overline C_k^i(q_n,t),i=1,2$ are surjective, for all $t\in[\frac{1}{10}s_n,s_n].$
\end{pro}

Before we give the proof, first we give a direct corollary of (2), which shows that the $\e$-process stated in the previous subsection has to stop at a finite step for any $k$.

\begin{cor}
 For any $k$ and $\e<\e_0$, the $\e$ process has to stop at a finite step.
\end{cor}

\nd Since $H^2(E_k\cap \overline B(0,1))<2\pi$, there exists $n_k>0$ such that
\be \inf_{q\in\R^4} H^2(P_k\bs D(q,s_{n_k}))>H^2(E_k).\ee
Then our process need to stop before the step $n_k$, because otherwise, we use the term (2) in Proposition \ref{reg}, for $t=s_n$, and get the disjoint decomposition
\be E_k=[E_k\cap D(q_{n_k},s_{n_k})]\cup G_{s_{n_k}}^1\cup G_{s_{n_k}}^2,\ee
therefore
\be \begin{split}
     H^2(E_k)&\ge H^2(G_{s_{n_k}}^1)+H^2(G_{s_{n_k}}^2)\ge H^2[p_k^1(G_{s_{n_k}}^1)]+H^2[p_k^2(G_{s_{n_k}}^2)]\\
&\ge H^2(P_k\bs D(q_{nk},s_{n_k})>H^2(E_k),
    \end{split}
\ee
which leads to a contradiction.\qed

Now we are going to prove Proposition 5.41.

\noindent Proof of Proposition 5.41.

For (1), notice that every topological minimal set is an Almgren minimal sets, hence (1) and (5.44) are direct corollaries of the proposition 6.1 (1). 

As a result of (1), we know that (5.45) in (2) is true if we replace all the $D_k(0,1)$ with $D_k(0,\frac{39}{40})$. Hence we have to prove that 
\be P_k^i\cap D(0,1)\bs D(0,\frac{39}{40})\subset p_k^i(G_t^i).\ee

We prove it for $i=1$ for example. The other case is the same. Denote by $P$ the plane orthogonal to $P_k^1$.

We know that $G_t^2$ is very closed to $P_k^2$, and when $k$ is large, $P_k^1$ and $P_k^2$ are almost orthogonal, hence the projection of $P_k^2\cap D(0,1)$ under $p_k^1$ is far away from $P_k^1\cap D(0,1)\bs D(0,\frac{39}{40}).$ On the other hand, the projection of the part $E_k\cap D(q_n,t)$ under $p_k^1$ is always contained in  $D(q_n,t)$, hence also far away from $P_k^1\cap D(0,1)\bs D(0,\frac{39}{40}).$ So (5.50) is equivalent to say that
\be P_k^1\cap D(0,1)\bs D(0,\frac{39}{40})\subset p_k^1(E_k\cap\overline D(0,1)).\ee

We are going to prove a stronger one, that is
\be P_k^1\cap \overline B(0,1)\subset p_k^1(E_k\cap\overline B(0,1)).\ee

So suppose that (5.52) is not true. That is, there exists $x\in P_k^1\cap \overline B(0,1)$ such that ${p_k^1}^{-1}\{x\}\cap (E_k\cap\overline B(0,1))=\emptyset.$ In other words, ${p_k^1}^{-1}\{x\}\cap\overline B(0,1)=(P+x)\cap\overline B(0,1)$ does not meet $E_k$. But $E_k$ is closed, so there exists $\d>0$ such that the neighborhood $B((P+x)\cap\overline B(0,1),\d)$ does not meet $E_k$. Denote by $S=(P+x)\cap\partial B(0,1+\d)$, then $S$ is a circle that does not meet $E_k$, and $S$ is zero in $H_1(\R^4\bs E_k)$, because it is the boundary of the disc $(P+x)\cap\partial B(0,1)\subset\R^4\bs E_k$.

However, $S$ is non zero in $H_1(\R^4\bs P_k^1)$, hence is non zero in $H_1(\R^4\bs P_k)$. This contradicts the fact that $E_k$ is a topological competitor of $P_k$. 

We have thus (5.52), which gives (5.51), and hence (5.50). 

So we get (2).

To prove (3), the idea is almost the same as above. We prove it for $i=1$ for example, denote still by $P$ the orthogonal plane of $P_k^1$. Suppose that there exists $x\in P_k^1\cap \overline C_k^1(q_n,t)$ satisfying ${p_k^1}^{-1}(x)\cap\overline D_k(q_n,t)\cap E_k=\emptyset$. Equivalently, 
\be (P+x)\cap\overline D_k(q_n,t)\cap E_k=\emptyset.\ee
Denote by $S=(P+x)\cap\partial D_k(q_n,t)$. Then $S$ is zero in $\R^4\bs E_k$, since it is the boundary of the disc $(P+x)\cap\overline D_k(q_n,t)$.
But since the $\e$-process does not stop at step $n$, outside $D_k(q_n,t)$, $E_k$ is composed of two disjoint pieces that are $\e$ closed to $P_k^1$ and $P_k^2$ respectively, and $d(q_n,0)\le 12\e$, so in fact we can deform our $S$ to any circle $S_y=(P+y)\cap \partial B(y,\frac 12)$ for all $y\in P_k^1\bs\overline B(0,1)$. Hence such a $S_y$ is zero in $H_1(\R^4\bs E_k)$. Again, we know that such a $S_y$ is non zero in $H_1(\R^4\bs P_k)$, which contradicts the fact that $E_k$ is a topological competitor of $P_k$.

Thus we get (3). And the proof of Proposition 5.41 is finished.\qed

\subsection{Conclusion}

The rest of the proof is exactly the same as the proof of the Almgren minimality of the union of two almost orthogonal planes in \cite{2p}. So we will only roughly describe what happens, without much detail.

So fix $k$ large and $\e$ small enough. Then outside $D_k(o_k,\frac{1}{10}r_k)$, $E_k$ is composed of two $C^1$ graphs $G^1,G^2$ on $P_k^1$ and $P_k^2$, and inside $D_k(o_k,r_k)$, $E_k$ is $\e r_k$ far from any translation of $P_k$. By a compactness argument we know that the part $E_k\cap D(o_k,r_k)\bs D_k(o_k,\frac{1}{10}r_k)$ is composed of two $C^1$ graphs that is far from any translation of $P_k$ (c.f. \cite{2p} Proposition 8.1). This makes one of the two graphs $G^1,G^2$, say $G^1$, oscillate of order $C(\e) r_k$, so that by an argument of harmonic extension (c.f. \cite{2p} Section 7 and 8) we know that the measure of $G^1$ is at least $C(\e)r_k^2$ more than its projection to $P_k^1$, that is, 
\be H^2(G^1)\ge H^2(D(o_k,\frac{1}{10}r_k)\cap P_k^1)+C(\e)r_k^2.\ee

For $G^2$, this is a graph on $D(o_k,\frac{1}{10}r_k)\cap P_k^2$, thus 
\be H^2(G^1)\ge H^2(D(o_k,\frac{1}{10}r_k)\cap P_k^2).\ee

For the part of $E_k$ inside $D(o_k,\frac{1}{10}r_k)$, Lemmas 5.11 and \ref{projection} gives
\be H^2(E_k\cap D(o_k,\frac{1}{10}r_k))\ge (1+2\cos\a_k)^{-1}H^2(P_k\cap D(o_k,\frac{1}{10}r_k)),\ee
where $\a_k$ is the first characteristic angle between $P_k^1$ and $P_k^2$. We sum over (5.54)-(5.56), and get
\be H^2(E_k)\ge H^2(P_k\cap D(0,1))+[C(\e)-C\cos\a_k]r_k^2,\ee
where $C(\e)$ depends only on $\e$. Thus when $k\to\infty$, $\cos\a_k\to 0$. Hence for $k$ large enough, (5.57) gives
\be H^2(E_k)>H^2(P_k\cap D(0,1))=2\pi,\ee
 which contradicts (5.25). So we get the desired contradiction, and this completes the proof of Theorem 5.1.\qed
 
 \begin{cor}Let $P^1$ and $P^2$ be two $m-$planes in $\R^{2m}$, with $m>2$. Suppose the $m$ characteristic angles between $P^1$ and $P^2$ are
 $\a_1=\a_2=\cdots=\a_{m-2}=0$ and $\theta<\a_{m-1}\le \a_m$, where $\theta$ is the $\theta$ in Theorem 5.1. Then their union $P^1\cup P^2$ is a $m-$dimensional topological minimal set.
 \end{cor}
 
 \nd This is the direct corollary of Theorem 5.1 and Proposition 3.23.\qed
 
 \begin{rem}We can extend the results of Theorem 5.1 and Corollary 5.59 to unions of $n$ $m-$dimensional planes. by the same idea. See \cite{XY10} for more detail.
 \end{rem}
 
 \begin{rem}Note that in \cite{XY10}, we have proved that the almost orthogonal union of two $m-$planes is Al-minimal. Now by Corollaries 3.17 and 5.59, we have found another family of unions of $m-$planes that are Al-minimal, which are non-transversal unions, and hence far from almost orthogonal. Intuitively we could say something like "interpolation" between these two cases, for example, unions of two $m-$planes with characteristic angles between these two cases are minimal. But up to now we do not know how to prove this.

 \end{rem}

\renewcommand\refname{References}
\bibliographystyle{plain}
\bibliography{reference}

\end{document}